\documentclass[oneside,notitlepage,12pt]{article}

\usepackage{amssymb}
\usepackage[leqno]{amsmath}
\usepackage{amsfonts}
\usepackage{amsopn}
\usepackage{amstext}
\usepackage{amsthm}

\usepackage{tikz}
\usetikzlibrary{cd}

\usepackage{enumitem}

\usepackage{verbatim}
\usepackage[colorlinks]{hyperref}
\usepackage{makeidx}

\newcommand{\define}[2]{{\em #1}\index{#2}}
\newcommand{\defined}[1]{\define{#1}{#1}}

\usepackage{calrsfs}
\usepackage{fourier-orns}
\usepackage{clock} %\ClockFrametrue\ClockStyle2
\usepackage[alpine, weather]{ifsym}

\providecommand{\cal}{\mathcal}
\renewcommand{\Bbb}{\mathbb}
\renewcommand{\frak}{\mathfrak}
\newenvironment{pf}{\begin{proof}}{\end{proof}}

\newcommand{\Dee}{{\cal{D}}}
\newcommand{\Ef}{{\cal{F}}}

\newcommand{\El}{{\cal{L}}}

\newcommand{\Yu}{{\cal{U}}}
\newcommand{\Vee}{{\cal{V}}}
\newcommand{\Wu}{{\cal{W}}}
\newcommand{\Zee}{{\Bbb{Z}}}
\newcommand{\Emm}{{\frak{M}}}

\newcommand{\Qyu}{{\Bbb{Q}}}
\newcommand{\Err}{{\Bbb{R}}}

\newcommand{\lam}{{\lambda}}

\newcommand{\sig}{\sigma}

\renewcommand{\phi}{\varphi}
\renewcommand{\rho}{\varrho}

\newcommand{\rest}{\restriction}

\newcommand{\ntr}{{n\in\omega}}

\newcommand{\loe}{\leq}
\newcommand{\goe}{\geq}

\newcommand{\rloe}{\preccurlyeq}

\newcommand{\subs}{\subseteq}
\newcommand{\sups}{\supseteq}
\newcommand{\nnempty}{\ne\emptyset}

\renewcommand{\iff}{\Longleftrightarrow}

\newcommand{\cl}{\operatorname{cl}}

\newcommand{\id}[1]{{\operatorname{i\!d}_{#1}}} % identity morphism
\newcommand{\cf}{\operatorname{cf}}
\newcommand{\dom}{\operatorname{dom}}
\newcommand{\cod}{\operatorname{cod}}

\newcommand{\oraz}{\qquad\text{and}\qquad}

\newcommand{\meet}{\wedge}

\newcommand{\Land}{\;\&\;}

\newcommand{\bd}{\partial}

\newcommand{\concat}{{}^\smallfrown}

\newtheorem{tw}{Theorem}[section]
\newtheorem{wn}[tw]{Corollary}
\newtheorem{lm}[tw]{Lemma}
\newtheorem{prop}[tw]{Proposition}
\newtheorem{claim}[tw]{Claim}
\theoremstyle{definition}
\newtheorem{df}[tw]{Definition}
\newtheorem{ex}[tw]{Example}

\theoremstyle{remark}

\newcommand{\setof}[2]{\{#1\colon #2\}}

\newcommand{\sett}[2]{\{#1\}_{#2}}
\newcommand{\sn}[1]{\{#1\}} % singleton
\newcommand{\dn}[2]{\{#1,#2\}} % doubleton
\newcommand{\pair}[2]{\langle #1, #2 \rangle} % pair
\newcommand{\triple}[3]{\langle #1, #2, #3 \rangle} % triple
\newcommand{\map}[3]{#1\colon #2 \to #3} % A function
\newcommand{\Cantor}{2^\omega}
\newcommand{\cantree}{2^{<\omega}}

\newcommand{\fra}{Fra\"iss\'e}

\providecommand{\nat}{\omega}

\newcommand{\ciag}[1]{{\sett{{#1}_n}{\ntr}}}

\newcommand{\ob}[1]{\operatorname{Obj}( #1 )}
\newcommand{\iso}{\approx}

\newcommand{\ciagi}[1]{\sig{#1}}

\newcommand{\finciagi}[1]{{#1}^{<\nat}}

\newcommand{\fK}{{\mathfrak{K}}}
\newcommand{\fL}{{\mathfrak{L}}}
\newcommand{\fM}{{\mathfrak{M}}}
\newcommand{\fC}{{\mathfrak{C}}}
\newcommand{\fS}{{\mathfrak{S}}}
\newcommand{\fT}{{\mathfrak{T}}}

\newcommand{\cmp}{\circ} % composition!!!
\newcommand{\sets}{\ensuremath{\mathfrak S\mathfrak e\mathfrak t}} % the category of sets
\newcommand{\BM}[1]{\operatorname{BM}\left(#1\right)}
\newcommand{\BMG}[2]{\BM {#1, #2} }

\newcommand{\below}[1]{{#1}^{\uparrow}}

\newcommand{\cont}{\ensuremath{\mathfrak c}}

\newcommand{\wek}[1]{{\vec{#1}}}
\providecommand{\ar}{\arrow}

\title{Weak \fra\ categories}
\author{
{\sc Wies{\l}aw Kubi\'s}\footnote{
Research supported by the GA\v{C}R project EXPRO  20-31529X and RVO: 67985840.
}
\\ \\
{\small Institute of Mathematics, Czech Academy of Sciences, Czechia}
}

\date{\clocktime\today}

\makeindex

\begin{document}

\maketitle

\begin{abstract}
We develop the theory of weak \fra\ categories, in which the crucial concept is the weak amalgamation property, discovered relatively recently in model theory. We show that, in a suitable framework, every weak \fra\ category has its unique {\em generic limit}, a special object in a bigger category, characterized by a certain variant of injectivity.
This significantly extends the present theory of \fra\ limits.

\ \\
\noindent
{MCS (2010):}
03C95, %(1980-now) Abstract model theory
18A30. %Limits and colimits (products, sums, directed limits, pushouts, fiber products, equalizers, kernels, ends and coends, etc.)

\noindent
{\it Keywords:} Weak amalgamation property, generic object, \fra\ limit.
\end{abstract}

\section*{Introduction}

Infinite mathematical structures are often built by constructing chains of finite (or in some sense \emph{finitary}) small \emph{blocks}.
For example, the random graph is typically viewed as the result of a random process of adding more and more vertices and edges, starting with the empty set. In the same manner, one can build a random metric space and any other first-order countably generated structure. Some of these structures have natural definitions or representations (like the random graph---introduced by Rado~\cite{Rado}, providing an explicit formula), while some other (like the Urysohn universal metric space~\cite{Urysohn}) can be explored only by looking at their properties or investigating suitable chains leading to the ``random" objects.
There is a lot of literature addressing random mathematical structures (including some dedicated journals). We would like to draw attention to another aspect of building structures from finitary blocks, addressing the question of existence of a \emph{generic} structure, that is, one that occurs most often.
This is strictly related to the theory of \emph{\fra\ limits}, namely, countable ultra-homogeneous structures.

One of our main tools is a category-theoretic variant of the classical Banach-Mazur game, introduced by Mazur around 1930.
The original Mazur's game was played with nonempty open intervals of the real line. The only rule was that in each step the interval must be contained in the previous one. The first player wins if the intersection of these intervals contains a point of a prescribed set of reals.
A slightly better and perhaps a bit more useful variant of this game is by exchanging the roles of the players.
This was done by Choquet (thus, the game is sometimes called \emph{Choquet game}) characterizing completeness of metric spaces. Actually, the full characterization of the existence of winning strategies is due to Oxtoby~\cite{Oxt}. We refer to Telg\'arsky~\cite{Rastislav} for a comprehensive survey of the Banach-Mazur game and its numerous variants.
When analyzing the ``essence" of this game, one immediately realizes that it can be played in an arbitrary category, where the inclusion of sets is replaced by abstract arrows. In particular, such a game was played many times in model theory, where the players alternately choose bigger and bigger models of the same type, see Hodges~\cite{HodBMG}.
We have already explored this game in~\cite{KubBMG} and \cite{KraKub}; the latter work characterizes the existence of a winning strategy for the second player, in the framework of model theory.
One of our goals here is to show that all these results can be proved by ``playing with arrows", namely, using the language and basic tools of category theory.
Our results extend the classical \fra\ theory, replacing the crucial amalgamation property by its weaker version.

We develop category-theoretic framework for the theory of generic limits of weak \fra\ classes. \fra\ theory belongs to the folklore of model theory, however actually it can be easily formulated in pure category theory.
The crucial point is the notion of \emph{amalgamation}, saying that two embeddings of a fixed object can be joined by further embeddings into a single one.
More precisely, for every two arrows $f,g$ with the same domain there should exist compatible arrows $f',g'$ with the same co-domain, such that $f' \cmp f = g' \cmp g$.
A significant relaxing of the amalgamation property, called the \emph{weak amalgamation property} was identified by Ivanov~\cite{Ivanov} and later independently by Kechris and Rosendal~\cite{KechRos} during their study of generic automorphisms in model theory.
Actually, this property is strictly related to the so-called \emph{pre-homogeneity}, which goes back to Pabion~\cite{Pab}.
It turns out that the weak amalgamation property is sufficient for constructing special unique objects satisfying certain variant of homogeneity.
We show how do it in pure category theory. 
We partially rely on the concepts and results of~\cite{Kub40}. One of our main goals is a general result on the existence and properties of special (called \emph{generic}) objects that are characterized up to isomorphism in terms of \emph{weak injectivity}. Some of our results extend~\cite{KraKub}, where the existence of generic objects in the classical model-theoretic setting has been characterized. Weak injectivity in the framework of pure category theory has been recently studied by Di Liberti~\cite{DiLi}, motivated by the results of~\cite{KraKub}.

\paragraph{}
This note is organized as follows.
The first sections contain rather technical results involving the weak amalgamation property, weak domination, and the crucial concept of a weak \fra\ sequence. Starting from Section~\ref{SectObjectsGenrkWeakInj}, we present an applicable framework involving the category $\fK$ of ``small" objects and the category $\fL$ of their limits. The crucial conditions (L0)--(L2) are formulated at the beginning of Section~\ref{SectObjectsGenrkWeakInj}. Our results can be easily applied to the model-theoretic setting, the only possible obstacle is that the property of ``being hereditary" is not seen from the category-theoretic perspective, however the weak amalgamation property allows ignoring this issue, at least to some extent.
This is discussed in Section~\ref{SectSlabeHomogenita} which also contains a short study of weak homogeneity. Section~\ref{SectBMGejmz} introduces the abstract Banach-Mazur game and explores its basic properties. The main result here is a characterization of weak injectivity, under the assumption that the base category is locally countable. We also show that the weak amalgamation property (which apparently is the main theme of this note) can be proved from the fact that the bigger category $\fL$ has a small enough cofinal subcategory.
Finally, the last section contains applications, including several relevant examples.

%\separator

\paragraph{Connections to the literature.}

The book of Hodges~\cite{HodMT} contains classical \fra\ theory for classes of finitely generated models (see the first section of Chapter 7 in~\cite{HodMT}).
The book of Ghilardi \& Zawadowski~\cite{GhiZaw} essentially uses the amalgamation property, defined in the context of monics. Furthermore, one can find in their Chapter 7 a result (Prop. 7.9) resembling the notion of a \fra\ sequence.
The paper of Kirby~\cite{Kirby} contains a strengthening of Droste \& G\"obel result~\cite{DroGoe}, significantly relaxing the cardinality assumption.
Finally, the work of Caramello~\cite{CaramO}
contains \fra's construction in the general setting involving domination (the notion introduced earlier by the present author), applying it in topos theory.

%\mrak{This needs to be expanded.}

\tableofcontents

\section{Preliminaries}

We start with some basic definitions and notation which will be needed in this note. For undefined notions concerning category theory we refer to Mac Lane~\cite{MacLane}.

Let $\fK$ be a category. The class of $\fK$-objects will be denoted by $\ob{\fK}$.
Given $a, b \in \ob{\fK}$, the set of all arrows from $a$ to $b$ will be denoted by $\fK(a,b)$.
The identity of a $\fK$-object $a$ will be denoted by $\id a$.
Slightly abusing notation and supporting the ideology that arrows are more important than objects, we will use the letter $\fK$ to denote the class of all $\fK$-arrows. In other words,
$$\fK = \bigcup_{a,b \in \ob{\fK}}\fK(a,b).$$
One of the axioms of a category says that $\fK(a,b) \cap \fK(a',b') = \emptyset$ whenever $\pair ab \ne \pair{a'}{b'}$.
Thus, given $f \in \fK$, there are uniquely determined objects $a, b$ such that $f \in \fK(a,b)$.
In this case $a$ is called the \define{domain}{domain} of $f$, denoted by $\dom(f)$, while $b$ is called the \define{co-domain}{co-domain} of $f$, denoted by $\cod(f)$.
The composition of arrows $f$ and $g$ will be denoted by $f \cmp g$. The composition makes sense if and only if $\dom(f) = \cod(g)$.
The fact that $f$ is an arrow with domain $x$ and co-domain $y$ will often be written as $\map f x y$.
Recall that $f$ is a \define{monic}{monic} (also called a \emph{monomorphism}) if for every compatible arrows $g_0, g_1$ with the same domain, the equation $f \cmp g_0 = f \cmp g_1$ implies $g_0 = g_1$. When we say that certain arrows are \define{compatible}{compatible arrows}, we mean that certain expressions involving composition of these arrows make sense. For instance, saying that ``$f$, $g$ are compatible" may mean that $f \cmp g$ makes sense or $g \cmp f$ makes sense, however this will always be clear from the context. In the definition of a monic, compatibility of $g_0$, $g_1$ means $\cod(g_0) = \dom(f) = \cod(g_1)$.

We shall use standard set-theoretic notation. In particular, as we have seen above, $\pair xy$ denotes an ordered pair. The letter $\nat$ denotes the set of all natural numbers (starting with zero) which at the same time is the first infinite ordinal and cardinal number, often denoted by $\aleph_0$.
Note that $\nat$ is also treated as a category, namely, $\pair n m$ is the unique arrow from $n$ to $m$ provided that $m \goe n$. More generally, every poset $\pair P \loe$ is a category, where the class of objects is $P$ and the class of arrows is $\loe$ (actually, it is enough to assume that the relation $\loe$ is reflexive and transitive, so that the axioms of a category are satisfied).
A covariant functor between posets is simply an \emph{increasing} (also called \emph{order preserving}) mapping.

Note that the set $\nat$ can also be treated as a monoid category (a category with one object), where the composition is addition and the identity arrow is $0$. In this note we shall never treat $\nat$ as a monoid.

By a \define{sequence}{sequence} in a category $\fK$ we mean a covariant functor from $\nat$ into $\fK$.
In order to make our notation economical (minimizing the number of symbols), a sequence will be denoted by 
$\map{\vec x}{\nat}{\fK}$ and, in turn, $x_n$ will denote the object $\vec x(n)$ and $x_n^m$ will denote the bonding arrow from $x_n$ to $x_m$ (formally, $x_n^m = \vec x (\pair n m)$, where $n \loe m$).
In other words, a sequence in $\fK$ will be encoded by a single letter with the vector symbol above (e.g., $\vec x$, $\vec a$, $\vec z$) and then suitable subscripts and superscripts will indicate its objects and bonding arrows.
Given a sequence $\vec x$, its colimit (possibly in a bigger category) will be denoted by $\pair{X}{\sett{x_n^\infty}{\ntr}}$ and we shall also write $X = \lim \vec x$.
Recall that $\sett{x_n^\infty}{\ntr}$ is the \defined{colimiting co-cone}, namely, it is a co-cone in the sense that $x_n^\infty = x_m^\infty \cmp x_n^m$ whenever $n < m$ and, given another co-cone $\sett{f_n}{\ntr}$ into a fixed object $Y$, there is a unique arrow $\map g X Y$ satisfying $g \cmp x_n^\infty = f_n$ for every $\ntr$. This is the formal definition of the colimit of a sequence.

Given sequences $\map{\vec x}{\nat}{\fK}$, $\map{\vec y}{\nat}{\fK}$, an arrow from $\vec x$ to $\vec y$ is typically a natural transformation. This is not good enough for our purposes.
Namely, we need to take into account more general arrows, specifically, natural transformations into subsequences. Thus, a \emph{transformation} from a sequence $\vec{x}$ to a sequence $\vec{y}$ is defined to be a natural transformation $\vec f$ from $\vec{x}$ to $\vec{y} \cmp \tilde{f}$, where $\map {\tilde{f}} \nat \nat$ is an increasing mapping. We denote by $f_n$ the arrow from $x_n$ to $y_{\tilde{f}(n)}$.
We have to identify transformations leading to the same colimit, namely, two transformations from $\vec{x}$ to $\vec{y}$ are \emph{equivalent} if the diagram consisting of both sequences and both transformations is commutative.
Intuitively, two transformations are equivalent if one can arrive at the same transformation by ``correcting" the corresponding increasing mappings of $\nat$.
For example, the identity of $\vec{x}$ is equivalent to the transformation sending $x_n$ to $x_{\psi(n)}$ via $x_n^{\psi(n)}$, where $\map \psi \nat \nat$ is any increasing mapping.
It is rather clear that this is an equivalence relation and identifying equivalent transformations we obtain a category structure on all sequences. Thus, from now on an \defined{arrow of sequences} will be an equivalence class of a transformation, as defined above.
The category of all sequences in $\fK$ with arrows defined above will be denoted by \defined{$\ciagi \fK$}.
Note that two sequences $\vec a$ and $\vec b$ are \define{isomorphic}{isomorphic sequences} if and only if there exist two transformations $\map{\vec g}{\vec a}{\vec b}$, $\map{\vec h}{\vec b}{\vec a}$ such that the diagram,
$$\begin{tikzcd}
	a_0 \ar[rr] \ar[dr, "f_0"'] && a_{\tilde{g}(\tilde{f}(0))} \ar[dr, "f_{\tilde{g}(\tilde{f}(0))}"'] \ar[rr] && a_{\tilde{g}(\tilde{f}(\tilde{g}(\tilde{f}(0)))} \ar[r] \ar[dr] & \cdots \\
	b_0 \ar[r] & b_{\tilde{f}(0)} \ar[rr] \ar[ur, "g_{\tilde{f}(0)}"] && b_{\tilde{f}(\tilde{g}(\tilde{f}(0))} \ar[ur, "g_{\tilde{f}(\tilde{g}(\tilde{f}(0))}"'] \ar[rr] && \cdots
\end{tikzcd}$$
consisting of both sequences and both transformations, is commutative. The equivalence classes of $\vec g$ and $\vec h$ provide a concrete isomorphism together with its inverse.

The concepts defined above are illustrated by the following simple example.

\begin{ex}
	Let $\fK$ be the category of finite sets with one-to-one mappings or, more generally, finite structures with embeddings. Particular types of sequences are chains $x_0 \subs x_1 \subs \cdots$, meaning that $x_n$ is a substructure of $x_{n+1}$ for each $\ntr$.
	Given two such chains $\vec x$, $\vec y$, we see that equivalent transformations from $\vec x$ to $\vec y$ are those leading to the same embedding of $X := \bigcup_{\ntr}x_n$ into $Y := \bigcup_{\ntr}y_n$. Note that every embedding of $X$ into $Y$ is the result of a transformation, as the image of every $x_n$ is contained in some $y_{\phi(n)}$ and $\phi$ may obviously be increasing.
	Finally, note that $X$ and $Y$ are colimits of the chains $\vec x$ and $\vec y$, respectively, in the category of all (or just countable) sets. Thus, $\ciagi{\fK}$ may be naturally identified with the category of countable structures that are unions of chains in $\fK$, where the arrows are arbitrary embeddings. Obviously, isomorphisms of chains correspond to bijections between their unions.
\end{ex}

We now recall one of the crucial concepts from abstract \fra\ theory that will be used later several times.

\begin{df}
	Following~\cite{Kub40}, we say that a subcategory $\fS \subs \fK$ is \define{dominating}{dominating subcategory} if the following conditions are satisfied.
	\begin{enumerate}[itemsep=0pt]
		\item[(C)] For every $x \in \ob{\fK}$ there is $f \in \fK$ such that $\dom(f) = x$ and $\cod(f) \in \ob{\fS}$.
		\item[(D)] For every $y \in \ob{\fS}$, for every $\fK$-arrow $\map f y z$ there is a $\fK$-arrow $\map g z u$ such that $g \cmp f \in \fS$ (in particular, $u \in \ob{\fS}$).
	\end{enumerate}
\end{df}

This definition makes sense for arbitrary families of arrows, however, any infinite dominating family generates a subcategory of the same cardinality and with the same objects, therefore domination is preserved.
The key point in \fra\ theory is finding countable dominating subcategories.
A subcategory $\fS$ satisfying condition (C) is called \define{cofinal}{cofinal subcategory}.

\section{Weak domination and weak amalgamations}

The amalgamation property is a well known concept in algebra and model theory. In category theory it is sometimes stated for monics only (see, e.g., \cite[p. 40]{GhiZaw}), although we do not see any formal reason for this, as all monics of a given category form a subcategory. On the other hand, in Sections~\ref{SectObjectsGenrkWeakInj}, \ref{SectSlabeHomogenita} we shall indeed assume that all arrows are monic.
Nevertheless, we need an important (and perhaps ultimate) weakening of the amalgamation property, therefore below we state somewhat technical variants of ``localized'' amalgamation properties.

\begin{df}
	Let $\fK$ be a fixed category.
	We shall say that $\fK$ has the \emph{amalgamation property at} $z \in \ob{\fK}$
	if for every $\fK$-arrows $\map f z x$, $\map g z y$ there exist $\fK$-arrows $\map{f'}{x}{w}$, $\map{g'}{y}{w}$ satisfying $f' \cmp f = g' \cmp g$. Such an object $z$ is also called an \emph{amalgamation base} (typically in model theory).
	Recall that $\fK$ has the \define{amalgamation property}{amalgamation property} (briefly: \emph{AP}) if it has the amalgamation property at every $z \in \ob{\fK}$.
\end{df}

A natural and important weakening is as follows.

\begin{figure}[b]
$$\begin{tikzcd}
& & y \ar[r, "g'"]& w \\
& & z' \ar[r,"f"'] \ar[u, "g"] & x \ar[u, "f'"'] \\
z \ar[rru, "e"] & & &
\end{tikzcd}$$
\caption{Cofinal amalgamation} \label{FigOneorJedn}	
\end{figure}
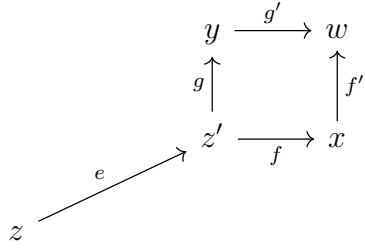

\begin{df}
	We say that $\fK$ has the \define{cofinal amalgamation property}{cofinal amalgamation property} (briefly: \emph{CAP}) if for every $z \in \ob{\fK}$ there exists a $\fK$-arrow $\map e z {z'}$ such that $\fK$ has the amalgamation property at $z'$ (see Fig.~\ref{FigOneorJedn}).
\end{df}

\begin{prop}
	A category has the cofinal amalgamation property if and only if it has a dominating subcategory with the amalgamation property.
\end{prop}

\begin{pf}
	Assume $\fK$ has the CAP and let $\fK_0$ be the full subcategory of $\fK$ such that
	$$\ob{\fK_0} = \setof{z \in \ob{\fK}}{\fK \text{ has the AP at } z}.$$
	We check that $\fK_0$ dominates $\fK$.
	CAP says that $\fK_0$ is cofinal in $\fK$, that is, (C) holds. As $\fK_0$ is full, (D) follows from (C).
	
	Now suppose that $\fS$ is a dominating subcategory of $\fK$ and $\fS$ has the AP.
	Fix $z \in \ob{\fK}$ and using (C) choose a $\fK$-arrow $\map e z u$ such that $u \in \ob{\fS}$.
	Fix $\fK$-arrows $\map f u x$, $\map g u y$.
	Using (D), find $\fK$-arrows $\map{f'}{x}{x'}$, $\map{g'}{y}{y'}$ such that $f' \cmp f, g' \cmp g \in \fS$.
	Applying the AP, find $\fS$-arrows $\map{f''}{x'}{w}$, $\map{g''}{y'}{w}$ such that $f'' \cmp f' \cmp f = g'' \cmp g' \cmp g$. This shows that $\fK$ has the AP at $u$.
\end{pf}

The above proposition shows that, from the category-theoretic point of view, cofinal AP is not much different from AP, as long as we agree to switch to a dominating subcategory.
Below is a significant and important weakening of the cofinal AP. In model theory, it was explicitly used first by Ivanov~\cite{Ivanov}, later by Kechris and Rosendal~\cite{KechRos}, and recently by Krawczyk and the author~\cite{KraKub}.

\begin{df}
	Let $\fK$ be a category. We say that $\fK$ has the \define{weak amalgamation property}{weak amalgamation property} (briefly: WAP)\footnote{Ivanov~\cite{Ivanov} calls it the \emph{almost amalgamation property}, while we follow Kechris and Rosendal~\cite{KechRos} who use the adjective \emph{weak} instead of \emph{almost}. Actually, we have already considered the concept of \emph{almost amalgamations} in metric-enriched categories~\cite{Kub41}, where the meaning of ``almost'' is, roughly speaking, ``commuting with a small error''.} if for every $z\in \ob{\fK}$ there exists a $\fK$-arrow $\map e z{z'}$ such that for every $\fK$-arrows $\map f{z'}x$, $\map g{z'}y$ there are $\fK$-arrows $\map{f'}{x}{w}$, $\map{g'}{y}{w}$ satisfying
	$$f' \cmp f \cmp e = g' \cmp g \cmp e.$$
\end{df}

In other words, the square in the diagram shown in Fig.~\ref{FigOneorJedn} may not be commutative.
The arrow $e$ above will be called \define{amalgamable}{amalgamable arrow} in $\fK$. Thus, $\fK$ has the WAP if for every $\fK$-object $z$ there exists an amalgamable $\fK$-arrow with domain $z$.
Note also that saying ``$\fK$ has the AP at $z \in \ob{\fK}$" is precisely the same as saying ``$\id z$ is amalgamable in $\fK$".

\begin{lm}\label{Lmebriesdv}
	Let $e \in \fK$ be an amalgamable arrow. Then $i \cmp e$ and $e \cmp j$ are amalgamable for every compatible arrows $i, j \in \fK$.
\end{lm}

\begin{pf}
	Assume $\map e z u$ is amalgamable, $\map i u v$, and fix $\map f v x$, $\map g v y$. Then $\map{f \cmp i}{u}{x}$ and $\map{g \cmp i}{u}{y}$ and therefore there are $\map{f'}{x}{w}$ and $\map{g'}{y}{w}$ such that $f' \cmp f \cmp i \cmp e = g' \cmp g \cmp i \cmp e$. This shows that $i\cmp e$ is amalgamable.
	It is clear that $e \cmp j$ is amalgamable as long as $e$ is.
\end{pf}

As it happens, the weak amalgamation property, contrary to its stronger variants, is very stable. In order to state it precisely, we need the concept of weak domination.

\begin{df}
	We say that $\fS \subs \fK$ is \define{weakly dominating}{weakly dominating subcategory} if it is cofinal (namely, satisfies (C) from the definition of domination) and
	\begin{enumerate}
		\item[(W)] For every $y \in \ob{\fS}$ there exists $\map j y {y'}$ in $\fS$ such that for every $\fK$-arrow $\map f{y'}z$ there is a $\fK$-arrow $\map g z u$ satisfying $g \cmp f \cmp j \in \fS$.
	\end{enumerate}
\end{df}

Note that for a full subcategory the stronger condition (D) follows from (C), therefore in that case being weakly dominating is the same as being dominating.
The somewhat technical concept of weak domination should be more clear after 
looking at the next result. This is an important characterization of the weak amalgamation property, clearly showing its strong stability with respect to suitable subcategories.

\begin{prop}\label{PROPdwapunktsix}
	Let $\fK$ be a category. The following properties are equivalent.
	\begin{enumerate}[itemsep=0pt]
		\item[{\rm (a)}] $\fK$ has the weak amalgamation property.
		\item[{\rm (b)}] Every cofinal full subcategory of $\fK$ has the weak amalgamation property.
		\item[{\rm (c)}] $\fK$ has a cofinal full subcategory with the weak amalgamation property.
		\item[{\rm (d)}] $\fK$ is dominated by a subcategory with the weak amalgamation property.
		\item[{\rm (e)}] $\fK$ is weakly dominated by a subcategory with the weak amalgamation property.
	\end{enumerate}	
\end{prop}

\begin{pf}
	(a)$\implies$(b) Let $\fS$ be cofinal and full in $\fK$, and fix $z \in \ob{\fS}$.
	Find an amalgamable $\fK$-arrow $\map e z v$.
	Using domination, we may find a $\fK$-arrow $\map i v u$ such that $i \cmp e \in \fS$.
	By Lemma~\ref{Lmebriesdv}, $i \cmp e$ is amalgamable in $\fK$. We need to show that it is amalgamable in $\fS$.
	For this aim, fix $\fS$-arrows $\map f v x$, $\map g v y$. Applying the WAP, we find $\fK$-arrows $\map{f'}{x}{w}$, $\map{g'}{y}{w}$ such that $f' \cmp f \cmp i \cmp e = g' \cmp g \cmp i \cmp e$.
	Finally, using domination again, find a $\fK$-arrow $\map j w {w'}$ such that $w' \in \ob{\fS}$. Then $j \cmp f'$ and $j \cmp g'$ are $\fS$-arrows, because $\fS$ is a full subcategory of $\fK$ (this is the only place where we use fullness).
	Finally, we have
	$$(j \cmp f') \cmp f \cmp (i \cmp e) = (j \cmp g') \cmp g \cmp (i \cmp e).$$
		
	(b)$\implies$(c)$\implies$(d)$\implies$(e) Obvious.
	
	(e)$\implies$(a) Let $\fS$ be weakly dominating in $\fK$ and assume $\fS$ has the WAP.
	Fix $z \in \ob{\fK}$. First, find a $\fK$-arrow $\map i z u$ with $u \in \ob{\fS}$.
	Now find an $\fS$-arrow $\map e u v$ that is amalgamable in $\fS$.
	Let $\map j v {v'}$ be an $\fS$ satisfying the assertion of (W).
	By Lemma~\ref{Lmebriesdv}, it suffices to show that $j \cmp e$ is amalgamable in $\fK$.
	Fix $\fK$-arrows $\map f v x$, $\map g v y$.
	Using domination, find $\fK$-arrows $\map{f'}{x}{x'}$, $\map{g'}{y}{y'}$ such that $f' \cmp f \cmp j \in \fS$ and $g' \cmp g \cmp j \in \fS$.
	Using the fact that $e$ is amalgamable in $\fS$, we find $\fS$-arrows $\map{f''}{x'}{w}$ and $\map{g''}{y'}{w}$ satisfying 
	$$f'' \cmp (f' \cmp f \cmp j) \cmp e = g'' \cmp (g' \cmp g \cmp j) \cmp e.$$
	Thus $j \cmp e$ is amalgamable in $\fK$.
\end{pf}

From the results above we see that the WAP comes from the CAP by shifting the amalgamation property from the objects to the arrows. In practice, at least to the author's knowledge, there are not many natural examples of categories with WAP that fail the CAP. Below is perhaps the first such example, due Pouzet, contained in Pabion~\cite{Pab}.

\begin{ex}[Pouzet]\label{EXMauriceP}
	Let $\El$ be the class of all finite linearly ordered sets. Given $X \in \El$, define a ternary relation $R_X$ by
	$$R_X(x,y,z) \iff x < z \Land y < z \Land x \ne y.$$
	We now forget the linear orderings, replacing them by the relations defined above.
	Namely, let $\fK$ be the category whose objects are all structures $\pair{X}{R_X}$, where $X \in \El$ and arrows are the embeddings.
	Note that, given $\pair{X}{R_X}$, one can ``almost'' reconstruct the linear ordering of $X$, except the first two elements.
	It is easy to verify that $\fK$ has the weak amalgamation property, while it definitely fails the cofinal one, as one can extend any $\pair Z{R_Z}$ in two incompatible ways by adding a new element below the first two elements of $Z$.
\end{ex}

One can argue whether the example above is natural or not. We shall see later that it leads to the set of rational numbers $\Qyu$ described in a different language. Note that one can fully reconstruct the ordering from the relation $R_\Qyu$. From the model-theoretic point of view, the structures $\pair \Qyu <$ and $\pair{\Qyu}{R_\Qyu}$ are inter-definable (each one can be defined from the other).

Other model-theoretic examples distinguishing WAP from CAP can be found in~\cite{KKKP}. Recall that the WAP was identified in~\cite{Ivanov} and \cite{KechRos} as the crucial ingredient for characterizing the existence of generic automorphisms, although all the examples in these works satisfy the CAP.
Finally, let us mention that it is formally quite easy to ``kill'' the cofinal amalgamation property: 

\begin{ex}
	Assume $\fK$ is a category with the CAP but not the AP, moreover, for every $z \in \ob{\fK}$ there exists a $\fK$-arrow $\map{e}{z}{z'}$ such that $\fK$ fails the AP at $z'$.
	Now let $\fK'$ be the full subcategory of $\fK$ obtained by removing all objects $z$ such that $\fK$ has amalgamations at $z$.
	Then $\fK'$ has the weak amalgamation property (by Proposition~\ref{PROPdwapunktsix}) while it evidently fails the cofinal amalgamation property.
	
	A very concrete example is the category $\fK$ of all finite cycle-free graphs with embeddings.
	Note that $\fK$ has the amalgamation property at $z$ if and only if $z$ is connected. Indeed, if $z$ has at least two components then one can consider two inclusions $\map f z x$, $\map g z y$ such that $x$ comes from $z$ by adding a path of length two joining two fixed components $z_0$, $z_1$ of $z$ and $y$ comes from $z$ by adding a path of length three, joining $z_0$ and $z_1$ at the same vertices. By this way, any amalgamation of $f$ and $g$ contains a cycle.
	Thus, the category $\fK'$ of all disconnected cycle-free finite graphs with embeddings satisfies the WAP and fails the CAP.	
\end{ex}

The example above is not particularly natural, however it clearly exhibits the fact that the WAP is much more stable than the CAP.

\section{Weak \fra\ sequences}\label{SectDSGdogugo}

We now define the crucial concept of this note.

\begin{df}
	Let $\fK$ be a fixed category. A sequence $\map{\vec u}{\nat}{\fK}$ will be called a \define{weak \fra\ sequence}{weak \fra\ sequence} if the following conditions are satisfied.
	\begin{enumerate}[itemsep=0pt]
		\item[(G1)] For every $x \in \ob{\fK}$ there is $n$ such that $\fK(x,u_n) \nnempty$. 
		\item[(G2)] For every $n \in \nat$ there exists $m \goe n$ such that for every $\fK$-arrow $\map f {u_m} y$ there are $k \goe m$ and a $\fK$-arrow $\map g y {u_k}$ satisfying $g \cmp f \cmp u_n^{m} = u_n^k$.
		$$\begin{tikzcd}
			\cdots \ar[r] & u_n \ar[rr, "u_n^m"] & & u_m \ar[dr, "f"'] \ar[rr, "u_m^k"] & & u_k \ar[r] & \cdots \\
			& & & & y \ar[ur, "g"']
		\end{tikzcd}$$
	\end{enumerate}
\end{df}

Condition (G1) says that the image of $\vec u$ is cofinal in $\fK$. Condition (G2) looks a bit technical, although it is actually strictly connected with the weak amalgamation property:

\begin{lm}\label{LmLMwapdirs}
	Every category with a weak \fra\ sequence is directed and has the weak amalgamation property.
\end{lm}

\begin{pf}
	Let $\vec u$ be a weak \fra\ sequence in $\fK$.
	Condition (G1) clearly implies that $\fK$ is directed, as every two $\fK$-objects have arrows into a single $u_n$ for $n$ big enough.
	
	Fix $n \in \nat$ and let $m \goe n$ be as in (G2). We claim that $u_n^m$ is amalgamable in $\fK$.
	Indeed, if $\map{f_0}{u_m}{x_0}$ and $\map{f_1}{u_m}{x_1}$ are $\fK$-arrows then there are $k_0,k_1 \goe m$ and $\fK$-arrows $\map{g_0}{x_0}{u_{k_0}}$, $\map{g_1}{x_1}{u_{k_1}}$ such that $g_i \cmp f_i \cmp u_n^m = u_n^{k_i}$ for $i=0,1$.
	Let $k \goe \max\{k_0,k_1\}$.
	Then
	$$(u_{k_0}^k \cmp g_0) \cmp f_0 \cmp u_n^m = u_n^k = (u_{k_1}^k \cmp g_1) \cmp f_1 \cmp u_n^m.$$
	Now, if $z \in \ob{\fK}$ and $\map e z {u_n}$ is a $\fK$-arrow (which exists by (G1)), then $u_n^m \cmp e$ is amalgamable, by Lemma~\ref{Lmebriesdv}.
\end{pf}

\begin{lm}\label{LmCzteyrdwua}
	A category with a weak \fra\ sequence is weakly dominated by a countable subcategory, namely, the subcategory generated by the image of a weak \fra\ sequence.
\end{lm}

\begin{pf}
	Assume $\map{\vec u}{\nat}{\fK}$ is a weak \fra\ sequence in a category $\fK$.
	Let $\fS$ be the subcategory generated by the image of $\vec u$.
	By (G1), $\fS$ is cofinal in $\fK$.
	Fix $x \in \ob{\fK}$ and let $\map e x z$ be a $\fK$-arrow, where $z \in \ob{\fS}$. Then $z = u_n$ for some $n \in \nat$. Let $m > n$ be as in condition (G2). Then, by the proof of Lemma~\ref{LmLMwapdirs}, $u_n^m$ is amalgamable in $\fK$. Thus, condition (W) is satisfied.	
\end{pf}

Concerning the proof above, let us remark that if a sequence $\map{\vec u}{\nat}{\fK}$ is one-to-one on objects (i.e. $u_n \ne u_m$ whenever $n<m$) then its image is already a subcategory. On the other hand, it may actually happen that $u_n = u$ for every $\ntr$. Then the category generated by the image of $\vec u$ has just one object $u$ and its arrows are of the form
$$u_{n_1}^{m_1} \cmp \dots \cmp u_{n_k}^{m_k},$$
where $n_1 \loe m_1, \dots, n_k \loe m_k$, $k \in \nat$.

Perhaps the most extreme example here is the category $\sets$ of all sets with all possible mappings, where a weak \fra\ sequence is any sequence in which singletons (i.e., one-element sets) appear cofinally. A very concrete example is $\map{\vec u}{\nat}{\sets}$, where $u_n = \sn0$ and $u_n^m$ is the identity for every $n<m$.
As $\sn0$ is terminal in $\sets$, (G1) and (G2) are clearly satisfied.
Note that $\sn0$ with its identity is a trivial monoid dominating $\sets$.

\begin{lm}\label{LmCtiriTriss}
	Assume $\fS \subs \fK$ is weakly dominating and $\map{\vec u}{\nat}{\fS}$ is a weak \fra\ sequence in $\fS$. Then $\vec u$ is a weak \fra\ sequence in $\fK$.
\end{lm}

\begin{pf}
	It is clear that the image of $\vec u$ is cofinal in $\fK$.
	It remains to check (G2).
	
	Fix $n$ and let $m \goe n$ be such that (G2) holds in $\fS$, namely:
	\begin{enumerate}
		\item[(1)] For every $\fS$-arrow $\map f {u_m}y$ there are $k \goe m$ and an $\fS$-arrow $\map g y {u_k}$ such that $g \cmp f \cmp u_n^m = u_n^k$.
	\end{enumerate}
	Let $\map e {u_m} a$ be such that (W) holds, namely:
	\begin{enumerate}
		\item[(2)] For every $\fK$-arrow $\map f a x$ there is a $\fK$-arrow $\map g x y$ such that $g \cmp f \cmp e \in \fS$.
	\end{enumerate}
	Applying (1), find $k \goe m$ and $\map i a {u_k}$ such that $i \cmp e \cmp u_n^m = u_n^k$.
	
	Fix a $\fK$-arrow $\map f {u_k} x$. Then $\map{f \cmp i}{a}{x}$, therefore applying (2) we can find a $\fK$-arrow $\map g x y$ such that $h := g \cmp f \cmp i \cmp e \in \fS$.
	Applying (1) to the $\fS$-arrow $h$, we find $\ell \goe k$ and an $\fS$ arrow $\map j y {u_\ell}$ such that $j \cmp h \cmp u_n^m = u_n^\ell$.
	Finally, we have
	$$u_n^\ell = j \cmp h \cmp u_n^m = j \cmp g \cmp f \cmp i \cmp e \cmp u_n^m = (j \cmp g) \cmp f \cmp u_n^k,$$
	which shows (G2).
\end{pf}

The concept of a weak \fra\ sequence is (as the name suggests) a natural generalization of the notion of a \define{\fra\ sequence}{\fra\ sequence} from~\cite{Kub40}, where it is required that $m=n$ in condition (G2).
The existence of a \fra\ sequence obviously implies directedness and the cofinal amalgamation property. Furthermore, a sequence isomorphic to a \fra\ sequence may not be \fra\ (it remains to be weak \fra, as we shall see later).
On the other hand, we have:

\begin{prop}
	Assume $\fK$ has the amalgamation property and $\vec u$ is a weak \fra\ sequence in $\fK$. Then $\vec u$ is a \fra\ sequence in $\fK$.
\end{prop}

\begin{pf}
	Fix $n$ and a $\fK$-arrow $\map f {u_n} y$. Let $m \goe n$ be such that (G2) holds.
	Using the AP, find $\map {f'} y w$ and $\map g {u_m} w$ such that $f' \cmp f = g \cmp u_n^m$.
	Using (G2), find $\map h w {u_k}$ with $k \goe m$ such that $h \cmp g \cmp u_n^m = u_n^k$.
	Finally, $h \cmp f' \cmp f = u_n^k$.	
\end{pf}

As we have already mentioned, the property of being a \fra\ sequence is not stable under isomorphisms of sequences, unless the category in question has the AP. It turns out that the property of being weak \fra\ is very stable.

\begin{prop}
	Assume $\vec u$, $\vec v$ are isomorphic sequences in $\fK$. If $\vec u$ is weak \fra\ then so is $\vec v$.
\end{prop}

\begin{pf}
	Let $\map{\vec p}{\vec u}{\vec v}$ and $\map{\vec q}{\vec v}{\vec u}$ be arrows of sequences whose compositions are equivalent to the identities. Assume $\vec u$ is weak \fra.
	Obviously, $\vec v$ satisfies (G1). It remains to check that $\vec v$ satisfies (G2).
	
	Fix $n \in \nat$ and let $n' \goe n$ be such that $\map{q_n}{v_n}{u_{n'}}$. Let $m \goe n'$ be such that (G2) holds for $\vec u$, namely, for every $\map f {u_m} x$ there are $k \goe m$ and $\map g x {u_k}$ satisfying $g \cmp f \cmp u_{n'}^m = u_{n'}^k$.
	Let $m' \goe m$ be such that $\map{p_m}{u_m}{v_{m'}}$.
	Then $m' \goe n$. We claim that $m'$ is ``suitable" for condition (G2) concerning the sequence $\vec v$.
	For this aim, fix a $\fK$-arrow $\map f {v_{m'}} y$. Applying (G2) to the sequence $\vec u$ and to the arrow $f \cmp p_m$, we obtain $k \goe m$ and a $\fK$-arrow $\map g y {u_k}$ satisfying
	$$g \cmp f \cmp p_m \cmp u_{n'}^m = u_{n'}^k.$$
	Let $k' \goe k$ be such that $\map{p_k}{u_k}{v_{k'}}$.
	Note that
	$$p_m \cmp u_{n'}^m \cmp q_n = v_n^{m'} \oraz p_k \cmp u_{n'}^k \cmp q_n = v_n^{k'},$$ because the composition of $\vec p$ with $\vec q$ is equivalent to the identity of $\vec v$, as shown in the following diagram.
	$$\begin{tikzcd}
	\cdots \ar[rr] & & u_{n'} \ar[r]  & u_m \ar[rd, "p_m"] \ar[rr] & & u_k \ar[rr] \ar[rd, "p_k"] & & \cdots \\
	\cdots \ar[r] & v_n \ar[ru, "q_n"] \ar[rrr] & & & v_{m'} \ar[rr] \ar[rd, "f"'] & & v_{k'} \ar[r] & \cdots \\
	& & & & & y \ar[uu, crossing over, "g"', near start]
	\end{tikzcd}$$
	Thus
	$$(p_k \cmp g) \cmp f \cmp v_n^{m'} = p_k \cmp g \cmp f \cmp p_m \cmp u_{n'}^m \cmp q_n = p_k \cmp u_{n'}^k \cmp q_n = v_n^{k'}.$$
	This shows (G2) and completes the proof.
\end{pf}

We shall later see that two sequences which are weak \fra\ in the same category are necessarily isomorphic.
It remains to show their existence.
We shall say that $\fK$ is a \define{weak \fra\ category}{weak \fra\ category}, if it is directed, has the weak amalgamation property, and is weakly dominated by a countable subcategory.

\begin{tw}\label{THMczteryszest}
	Let $\fK$ be a category.
	The following properties are equivalent:
	\begin{enumerate}[itemsep=0pt]
		\item[{\rm(a)}] $\fK$ is a weak \fra\ category.
		\item[{\rm(b)}] There exists a weak \fra\ sequence in $\fK$.
	\end{enumerate}
\end{tw}

\begin{pf}
	Implication (b)$\implies$(a) is the content of Lemmas~\ref{LmLMwapdirs} and~\ref{LmCzteyrdwua}. It remains to show (a)$\implies$(b).
	By Lemma~\ref{LmCtiriTriss}, we may assume that $\fK$ itself is countable.
	In order to show the existence of a weak \fra\ sequence, we shall use the following simple claim, known in set theory as the Rasiowa-Sikorski Lemma:
	
	\begin{claim}
		Let $\pair P \loe$ be a partially ordered set and let $\Dee$ be a countable family of cofinal subsets of $P$.
		Then there exists a sequence $p_0 \loe p_1 \loe p_2 \loe \cdots$ in $P$ such that $D \cap \setof{p_n}{n \in \nat} \nnempty$ for every $D \in \Dee$.
	\end{claim}
	
	Let $\finciagi{\fK}$ denote the set of all finite sequences in $\fK$, that is, all covariant functors from $n = \{0,1, \dots, n-1\}$ into $\fK$, where $n \in \nat$ is arbitrary.
	We shall use the same convention as for infinite sequences, namely, if $\map{\vec x}{n}{\fK}$ then we shall write $x_i$ instead of $x(i)$ and $x_i^j$ instead of $x(\pair i j)$.
	Given $\vec a, \vec b \in \finciagi{\fK}$, define $\vec a \loe \vec b$ if $\vec b$ extends $\vec a$. Clearly, $\pair{\finciagi{\fK}}{\loe}$ is a partially ordered set. An increasing sequence in $\finciagi{\fK}$ gives rise to an infinite sequence in $\fK$, as long as it does not stabilize. Let $P$ be the subset of $\finciagi{\fK}$ consisting of all sequences $\map{\vec x}{n}{\fK}$ such that $x_i^j$ is amalgamable in $\fK$ whenever $i < j$.
	We shall work in the partially ordered set $\pair P \loe$.

	Given $x \in \ob{\fK}$, define $\Yu_x$ to be the set of all $\vec x \in P$ such that there is a $\fK$-arrow from $x$ to $x_i$ for some $i < \dom(\vec x)$.
	As $\fK$ is directed and has the weak AP, $\Yu_x$ is cofinal in $\pair{\finciagi{\fK}}{\loe}$.
	This follows from the fact that every $\fK$-arrow  can be prolonged to an amalgamable one (see Lemma~\ref{Lmebriesdv}).
	
	Fix $n \in \nat$ and $f \in \fK$. Define
	$\Vee_{n,f}$ to be the set of all $\vec x \in P$ such that $n+1 \in \dom(x)$ and the following implication holds:
	\begin{enumerate}
		\item[($*$)] If $x_{n+1} = \dom(f)$ then there are $k > n$ and $g \in \fK$ such that $g \cmp f \cmp x_{n}^{n+1} = x_n^m$.
	\end{enumerate}
	We check that $\Vee_{n,f}$ is cofinal in $P$.
	Fix $\vec a \in P$. First, we extend $\vec a$ by using amalgamable arrows so that $n+1 < \dom(\vec a)$.
	Now if $a_{n+1} \ne \dom(f)$ then already $\vec a \in \Vee_{n,f}$, so suppose $a_{n+1} = \dom(f)$.
	Let $k = \dom(\vec a)$ and assume $\map f {a_{n+1}} y$.
	Knowing that $a_n^{n+1}$ is amalgamable, we can find $\fK$-arrows $\map{g}{y}{w}$, $\map{h}{a_{k-1}}{w}$ such that
	$g \cmp f \cmp a_n^{n+1} = h \cmp a_{n+1}^{k-1} \cmp a_n^{n+1}$.
	Extend $\vec a$ by adding the arrow $h$ on the top, so that ($*$) holds. The extended sequence is a member of $\Vee_{n,f}$.
	This shows that $\Vee_{n,f}$ is cofinal in $\pair{P}{\loe}$.
	
	Finally, observe that a sequence $\vec{p_0} \loe \vec{p_1} \loe \vec{p_2} \loe \cdots$ satisfying the assertion of the Rasiowa-Sikorski Lemma (with $\Dee$ consisting of all possible $\Yu_x$ and $\Vee_{f,n}$) yields a weak \fra\ sequence in $\fK$.
	This completes the proof.
\end{pf}

\begin{df}
	A weak \fra\ sequence $\vec u$ is \define{normalized}{weak \fra\ sequence!-- normalized} if for every $n$ condition (G2) holds with $m = n+1$. More precisely, for every $n$, for every arrow $\map f {u_{n+1}} y$ there are $k>n$ and and an arrow $\map g y {u_k}$ such that $g \cmp f \cmp u_n^{n+1} = u_n^k$.
\end{df}
The sequence obtained in the proof above is normalized. Clearly, every weak \fra\ sequence contains a subsequence that is normalized.
In a normalized weak \fra\ sequence all non-identity  bonding arrows are amalgamable.
It turns out that the converse is true as well:

\begin{lm}
	Let $\vec u$ be a weak \fra\ sequence in $\fK$ such that $u_n^{n+1}$ is amalgamable for every $\ntr$.
	Then $\vec u$ is normalized.
\end{lm}

\begin{pf}
	Fix a $\fK$-arrow $\map f{u_{n+1}}y$. Let $m>n+1$ be as in condition (G2) applied to $n+1$ instead of $n$.
	Using the fact that $u_n^{n+1}$ is amalgamable, we find $\fK$-arrows $\map h y z$ and $\map{f'}{u_m}{z}$ such that $h \cmp f \cmp u_n^{n+1} = f' \cmp u_{n+1}^m \cmp u_n^{n+1}$.
	Using (G2), we find $k \goe m$ and a $\fK$-arrow $\map{g'}{z}{u_k}$ satisfying $g' \cmp f' \cmp u_{n+1}^m = u_{n+1}^k$.
	Let $g := g' \cmp h$. Then
	$$g \cmp f \cmp u_n^{n+1} = g' \cmp h \cmp f \cmp u_n^{n+1} = g' \cmp f' \cmp u_{n+1}^m \cmp u_n^{n+1} = u_{n+1}^k \cmp u_n^{n+1} = u_n^k,$$
	showing that $\vec u$ is normalized.
\end{pf}

The following fact will be essential for proving a variant of homogeneity of generic objects. The proof is a suitable adaptation of the back-and-forth argument.

\begin{lm}\label{LMverugbyrg}
	Assume $\vec u$, $\vec v$ are weak \fra\ sequences in $\fK$ such that $u_0^1$ is amalgamable and $\map{f}{u_1}{v_1}$ is a $\fK$-arrow.
	Then there exists an isomorphism of sequences $\map{\vec h}{\vec u}{\vec v}$ extending $f \cmp u_0^1$.
\end{lm}

\begin{pf}
	Passing to subsequences, we may assume $\vec u$, $\vec v$ are normalized and $u_0^1$ remains as it was, as it is already amalgamable.
	We construct the following (not necessarily commutative!) diagram
	$$\begin{tikzcd}
	u_{k_1} \ar[r] & u_{k_1+1} \ar[r] \ar[d, "f_1"'] & u_{k_2} \ar[r] & u_{k_2+1} \ar[d, "f_2"'] \ar[r] & u_{k_3} \ar[r] & \cdots \\
	v_0 \ar[r] & v_{\ell_1} \ar[r] & v_{\ell_1+1} \ar[u, "g_1"] \ar[r] & v_{\ell_2} \ar[r] & v_{\ell_2+1} \ar[u, "g_2"] \ar[r] & \cdots
	\end{tikzcd}$$
	in which $k_1 = 0$, $\ell_1 = 1$, and $f_1 = f$.
	Furthermore,
	\begin{enumerate}[itemsep=0pt]
		\item[(1)] $g_i \cmp v_{\ell_i}^{\ell_i+1} \cmp f_i \cmp u_{k_i}^{k_i+1} = u_{k_i}^{k_{i+1}}$,
		\item[(2)] $f_{j+1} \cmp u_{k_{j+1}}^{k_{j+1}+1} \cmp g_j \cmp v_{\ell_j}^{v_{\ell_j+1}} = v_{\ell_j}^{\ell_{j+1}}$
	\end{enumerate}
	holds for all $i,j \in \nat$.
	The construction is possible, because both sequences are normalized weak \fra, and hence (1), (2) are straightforward applications of the normalized variant of (G2).
	Define
	$$h_i = f_i \cmp u_{k_i}^{k_i+1} \oraz q_j = g_j \cmp v_{\ell_j}^{\ell_j+1}.$$
	Equations (1) and (2) give $q_i \cmp h_i = u_{k_i}^{k_{i+1}}$ and $h_{j+1} \cmp q_j = v_{\ell_j}^{\ell_{j+1}}$
	for $i,j \in \nat$.
	Thus $\vec h = \ciag{h}$ is an isomorphism from $\vec u$ to $\vec v$ and it extends $h_1 = f \cmp u_0^1$.
\end{pf}

\begin{wn}\label{WNCorUniquenesss}
	A category may have, up to isomorphism, at most one weak \fra\ sequence.
\end{wn}

\begin{pf}
	Let $\vec u$, $\vec v$ be weak \fra\ in $\fK$.
	Replacing them by subsequences, we may assume that they are normalized.
	By (G1), there exists a $\fK$-arrow $\map{f}{u_1}{v_k}$ for some $k$. Further refining $\vec v$, we may assume $k=1$. Now Lemma~\ref{LMverugbyrg} yields an isomorphism from $\vec u$ to $\vec v$.
\end{pf}

We finish this section by proving the following weakening of cofinality (in model theory usually called \emph{universality}).

\begin{lm}\label{LMsztyrydziesiync}
	Let $\vec u$ be a weak \fra\ sequence in $\fK$ and let $\vec x$ be a sequence in $\fK$ such that $x_n^{n+1}$ is amalgamable in $\fK$ for every $n \in \nat$.
	Then there exists a $\ciagi \fK$-arrow $\map{\vec e}{\vec x}{\vec u}$.
\end{lm}

\begin{pf}
	For simplicity, we assume that the sequence $\vec u$ is normalized.
	We construct inductively $\fK$-arrows $\map{e_n}{x_n}{u_{s(n)}}$ so that the following conditions are satisfied.
	\begin{enumerate}[itemsep=0pt]
		\item[(1)] $u_{s(n)}^{s(n+1)} \cmp e_n = e_{n+1} \cmp x_{n}^{n+1}$.
		\item[(2)] $e_n = e'_n \cmp x_n^{n+2}$ for some $\fK$-arrow $\map{e'_n}{x_{n+2}}{u_{s(n)}}$.
	\end{enumerate}
	We start with $e_0 = e'_0 \cmp x_0^2$, where $e'_0$ is an arbitrary $\fK$-arrow from $x_{n+2}$ into some $u_{s(0)}$, which exists by (G1).
	Suppose $e_0, \dots, e_n$ have been constructed.
	Let $\map{f}{x_{n+3}}{w}$ and $\map{g}{u_{s(n)+1}}{w}$ be $\fK$-arrows such that
	$$f \cmp x_{n+2}^{n+3} \cmp x_{n+1}^{n+2} = g \cmp u_{s(n)}^{s(n)+1} \cmp e'_n \cmp x_{n+1}^{n+2}.$$ This is possible, because $x_{n+1}^{n+2}$ is amalgamable in $\fK$.
	Using (G2) and the fact that $\vec u$ is normalized, we find a $\fK$-arrow $\map{h}{w}{u_{s(n+1)}}$, with $s(n+1) > s(n)$, such that $h \cmp g \cmp u_{s(n)}^{s(n)+1} = u_{s(n)}^{s(n+1)}$.
	Define $e'_{n+1} := h \cmp f$ and $e_{n+1} := e'_{n+1} \cmp x_{n+1}^{n+3}$.
	Then 
	$$e_{n+1} \cmp x_n^{n+1} = h \cmp f \cmp x_{n+1}^{n+3} \cmp x_n^{n+1} = h \cmp g \cmp u_{s(n)}^{s(n)+1} \cmp e'_n \cmp x_{n+1}^{n+2} \cmp x_n^{n+1} = u_{s(n)}^{s(n+1)} \cmp e_n.$$
	It follows that the construction can be carried out, obtaining a $\ciagi \fK$-arrow $\map{\vec e}{\vec x}{\vec u}$ with $\vec e = \ciag{e}$.	
\end{pf}

The results above show the importance of amalgamable arrows. One can say that the weak \fra\ theory comes from the usual one by moving the relevant concepts from the objects to the arrows. By this way we obtain a framework that is both more general and more robust, in the sense that it does not affect dominating subcategories.
Let us admit that in the category-theoretic approach to \fra\ theory, the cofinal amalgamation property plays the crucial role, as one can always restrict to a full cofinal subcategory. On the other hand, a full cofinal subcategory may fail the cofinal AP, as the example below shows. The weak version is much more stable, due to Proposition~\ref{PROPdwapunktsix}.

\begin{ex}\label{EXdisconsBushes}
	Let $\fT$ be the category of all finite cycle-free graphs (simple undirected graphs with no cycles) with embeddings. Then $\fT$ has the cofinal amalgamation property, namely it has amalgamations at $A \in \ob{\fT}$ if and only if $A$ is connected.
	Indeed, if $A$ is disconnected and $x,y \in A$ come from different components then adding a path joining $x$ and $y$ provides a cycle-free graph containing $A$; by adding two paths of different length we obtain two embeddings of $A$ that cannot be amalgamated.
	On the other hand, if $A$ is connected then every two embeddings of $A$ into cycle-free graphs can be amalgamated in the ``minimal'' way, namely not adding any unnecessary edge. More precisely, if $B_0, B_1$ are cycle-free, $A = B_0 \cap B_1$ then $B_0 \cup B_1$ is cycle-free as long as we do not add any edges between $B_0 \setminus A$ and $B_1 \setminus A$.
	
	Now let $\fK$ be the full subcategory of $\fT$ whose objects are precisely the disconnected graphs. By the arguments above, $\fK$ totally fails the cofinal AP. On the other hand, by Proposition~\ref{PROPdwapunktsix}, it still has the weak AP.
	
 \end{ex}

\section{Weakly injective objects}\label{SectObjectsGenrkWeakInj}

The previous section was somewhat technical, as we were working in the rather abstract category of sequences. We now prepare the setup suitable for exploring generic objects. For obvious reasons, they could be called \emph{generic limits} of weak \fra\ categories. Another possible and tempting name would be \emph{weak \fra\ limit}, however, in our opinion this would be a little bit inappropriate, because we only relax \fra's axioms, showing that the generic \emph{limit} is still unique and may only have weaker properties. After all, a weak \fra\ category may contain a weakly dominating \fra\ subcategory, having the same generic limit (see, e.g.,  Example~\ref{EXdisconsBushes} above). In any case, we shall avoid the word \emph{limit}, adapting the terminology from set-theoretic forcing, at some point calling the limit of a weak \fra\ sequence a \emph{generic object} (see Section~\ref{SectBMGejmz}). In this section we characterize those objects by a variant of injectivity.

As before, $\fK$ will denote a fixed category. Now we also assume that $\fL \sups \fK$ is a bigger category such that $\fK$ is full in $\fL$ and the following conditions are satisfied:
\index{(L0)--(L2)}
\begin{enumerate}[itemsep=0pt]
	\item[(L0)] All $\fL$-arrows are monic.
	\item[(L1)] Every sequence in $\fK$ has a colimit in $\fL$ and every $\fL$-object is the colimit of some sequence in $\fK$.
	\item[(L2)] Every $\fK$-object is \define{small}{small object} in the following sense: If $Y = \lim \vec y$, where $\vec y$ is a sequence in $\fK$, then for every $\fL$-arrow $\map f x Y$ there are $n$ and an $\fL$-arrow $\map {f'} x {y_n}$ such that $f = y_n^\infty \cmp f'$, where $y_n^\infty$ denotes the $n$th arrow from the colimiting co-cone.
\end{enumerate}
Concerning (L2), recall that there is a well-established concept of a \emph{finitely presented object} (sometimes called a \emph{compact object}), with a similar definition, using arbitrary functors from directed posets or even filtered categories (see e.g.~\cite{AdaRos}). On the other hand, condition (L2) describes exactly what we need, and nothing more. So, if every $\fK$-object is finitely presented in $\fL$ then (L2) holds, however the converse is not true: Consider the poset category $\fL$ consisting of all ordinals $\loe \omega_1$ and let $\fK$ consist of all ordinals of cofinality $\ne \omega$; then $\omega_1$ is not finitely presented in $\fL$, however it is small in the sense of (L2).

We shall use the following convention: The $\fL$-objects and $\fL$-arrows will be denoted by capital letters, while the $\fK$-objects and arrows will be denoted by small letters.

Typical examples of pairs $\pair \fK \fL$ satisfying (L0)--(L2) come from model theory: $\fK$ could be any class of finite (or just finitely generated) structures of a fixed first order language while $\fL$ should be the class of all structures isomorphic to the unions of countable chains of $\fK$-objects. The arrows in both categories are typically all embeddings.

It turns out that for every category $\fK$ in which all arrows are monic, the sequence category $\ciagi \fK$ can play the role of $\fL$, however in the applications one usually has in mind a more concrete and natural category satisfying (L0)--(L2). This is evident in Section~\ref{SectAppsyJ}, where we discuss sample applications.

\begin{df}
	We say that $U \in \ob{\fL}$ is \define{weakly $\fK$-injective}{weakly injective object} if 
	\begin{enumerate}[itemsep=0pt]
		\item[(U)] Every $\fK$-object has an $\fL$-arrow into $U$ (in other words: $\fL(x,U) \nnempty$ for every $x \in \ob{\fK}$).
		\item[(WI)] For every $\fL$-arrow $\map e a U$ there exists a $\fK$-arrow $\map i a b$ such that for every $\fK$-arrow $\map f b y$ there is an $\fL$-arrow $\map g y U$ satisfying $g \cmp f \cmp i = e$, as shown in the following diagram.
		$$\begin{tikzcd}
			a \ar[r, "i"] \ar[d, "e"'] & b \ar[r, "f"] & y \ar[dll, dashrightarrow, "g"] \\
			U
		\end{tikzcd}$$
	\end{enumerate}
\end{df}

As one can expect, this concept is strictly related to weak \fra\ sequences.
Recall that in this section we assume (L0)--(L2).

\begin{tw}\label{ThmMgenerikgnrc}
	Let $U = \lim \vec u$, where $\vec u$ is a sequence in $\fK$.
	Then $U$ is weakly $\fK$-injective if and only if $\vec u$ is a weak \fra\ sequence in $\fK$.
\end{tw}

\begin{pf}
	Assume first that $U$ is weakly $\fK$-injective.
	Condition (U) combined with (L2) shows that the sequence $\vec u$ satisfies (G1).
	In order to check (G2), fix $n \in \nat$ and apply the weak $\fK$-injectivity of $U$ to the arrow $\map{u_n^\infty}{u_n}{U}$.
	We obtain a $\fK$-arrow $\map i {u_n} b$ such that for every $\fL$-arrow $\map f b y$ there is an $\fL$-arrow $\map g y U$ satisfying $g \cmp f \cmp i = u_n^\infty$.
	Taking $f = \id b$, we obtain an $\fL$-arrow $\map j b U$ such that
	$$j \cmp i = u_n^\infty.$$
	Applying (L2), we get $m >n$ and a $\fK$-arrow $\map k b {u_m}$ such that $j = u_m^\infty \cmp k$.
	Thus
	$$u_m^\infty \cmp u_n^m = u_n^\infty = j \cmp i = u_m^\infty \cmp k \cmp i.$$
	By (L0), $u_m^\infty$ is a monic, therefore
	$$k \cmp i = u_n^m.$$
	We claim that $m$ is a witness for (G2).
	Fix a $\fK$-arrow $\map f{u_m}y$.
	Applying weak $\fK$-injectivity to the arrow $f \cmp k$, we find $\map g y U$ such that
	$$g \cmp f \cmp k \cmp i = u_n^\infty.$$
	Using (L2), we find $\ell>m$ and an $\fL$-arrow $\map {g'}y{u_\ell}$ such that $g = u_\ell^\infty \cmp g'$. Now we have
	$$u_\ell^\infty \cmp g' \cmp f \cmp u_n^m = u_\ell^\infty \cmp g' \cmp f \cmp k \cmp i = g \cmp f \cmp k \cmp i = u_n^\infty = u_\ell^\infty \cmp u_n^\ell.$$
	As $u_\ell^\infty$ is a monic, we conclude that $g' \cmp f \cmp u_n^m = u_n^\ell$, showing (G2).
	
	Now suppose that $\vec u$ is a weak \fra\ sequence in $\fK$.
	Then (G1) implies that $U$ satisfies (U).
	It remains to show that $U$ is weakly $\fK$-injective.
	Fix $\map e a U$. Using (L2), find $n$ and a $\fK$-arrow $\map {e'}a{u_n}$ such that $e = u_n^\infty \cmp e'$. Let $m>n$ be such that the assertion of (G2) holds.
	Define $i = u_n^m \cmp e'$.
	Fix a $\fK$-arrow $\map{f}{u_m}{y}$.
	There are $k\goe m$ and a $\fK$-arrow $\map {g'} y {u_k}$ such that $g' \cmp f \cmp u_n^m = u_n^k$. Let $g = u_k^\infty \cmp g'$.
	Then
	$$g \cmp f \cmp i = u_k^\infty \cmp g' \cmp f \cmp u_n^m \cmp e' = u_k^\infty \cmp u_n^k \cmp e' = u_n^\infty \cmp e' = e.$$
	Thus, $i$ witnesses the weak $\fK$-injectivity of $U$.
\end{pf}

Recall that $\fK$ has a weak \fra\ sequence if and only if it is a weak \fra\ category, i.e., it is directed, has the weak amalgamation property, and is weakly dominated by a countable subcategory.

\begin{wn}\label{WNCorpietdwwwaa}
	A weakly $\fK$-injective object exists if and only if $\fK$ is a weak \fra\ category.
\end{wn}

\begin{pf}
	If $\fK$ is a weak \fra\ category then it has a weak \fra\ sequence, whose colimit in $\fL$ is a weakly $\fK$-injective object by Theorem~\ref{ThmMgenerikgnrc}.
	Conversely, if $U$ is weakly $\fK$-injective then, by (L1), $U = \lim \vec u$ for some sequence $\vec u$ in $\fK$. By Theorem~\ref{ThmMgenerikgnrc}, the sequence $\vec u$ is weak \fra\ in $\fK$. Finally, by Theorem~\ref{THMczteryszest}, $\fK$ is a weak \fra\ category.
\end{pf}

\begin{wn}
	A weakly $\fK$-injective object, if exists, is unique up to isomorphism.	
\end{wn}

\begin{pf}
	Suppose $U$, $V$ are weakly $\fK$-injective. By (L1), $U = \lim \vec u$, $V = \lim \vec v$, where $\vec u$, $\vec v$ are sequences in $\fK$.
	By Theorem~\ref{ThmMgenerikgnrc}, both $\vec u$ and $\vec v$ are weak \fra\ in $\fK$.
	By Corollary~\ref{WNCorUniquenesss}, there exists an isomorphism from $\vec u$ to $\vec v$ in the category of sequences. This leads to an isomorphism between $U$ and $V$.
\end{pf}

\begin{wn}
	Let $U$ be a weakly $\fK$-injective object.
	If $X = \lim \vec x$, where $\vec x$ is a sequence in $\fK$ such that each bonding arrow $x_n^{n+1}$ is amalgamable in $\fK$, then there exists an $\fL$-arrow from $X$ to $U$.
\end{wn}

\begin{pf}
	Knowing that $U = \lim \vec u$, where $\vec u$ is a weak \fra\ sequence in $\fK$, it suffices to apply Lemma~\ref{LMsztyrydziesiync}.
\end{pf}

We now turn to the question of homogeneity.

\begin{tw}\label{THMigiuweg}
	Let $U$ be a weakly $\fK$-injective object and let $\map e a b$ be an amalgamable arrow in $\fK$.
	Then for every $\fL$-arrows $\map i b U$, $\map j b U$ there exists an automorphism $\map h U U$ satisfying $h \cmp i \cmp e = j \cmp e$.
\end{tw}

This is illustrated in the following diagram in which the triangle is not necessarily commutative.
$$
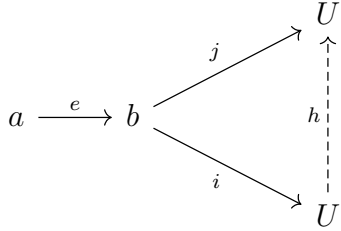
\begin{tikzcd}
& & & U \\
a \ar[r, "e"] & b \ar[rru, "j"] \ar[rrd, "i"'] & & \\
& & & U \ar[uu, dashrightarrow, "h"]
\end{tikzcd}$$

\begin{pf}
	Assume $U = \lim \vec u$, where $\vec u$ is a normalized weak \fra\ sequence in $\fK$.
	By (L2), there are $k,\ell \in \nat$ such that $i = u_k^\infty \cmp i'$ and $j = u_\ell^\infty \cmp j'$. We may assume that $\ell>0$, replacing $j'$ by $u_\ell^{\ell+1} \cmp j'$, if necessary. We may also assume that $i'$ is amalgamable, replacing it by $u_k^{k+1} \cmp i'$ (and increasing $k$), if necessary.
	Now
	$$\begin{tikzcd}
		a \ar[r, "e"] & b \ar[r, "i'"] & u_k \ar[r, "u_k^{k+1}"] & u_{k+1} \ar[r] & \cdots
	\end{tikzcd}$$
	and
	$$\begin{tikzcd}
	u_{\ell-1} \ar[r, "u_{\ell-1}^\ell"] & u_\ell \ar[r, "u_\ell^{\ell+1}"] & u_{\ell+1} \ar[r] & \cdots
	\end{tikzcd}$$
	are normalized weak \fra\ sequences and $\map{j'}{b}{u_\ell}$ is a $\fK$-arrow. By Lemma~\ref{LMverugbyrg}, there is an isomorphism of sequences $\vec h$ extending $j' \cmp e$. This leads to an isomorphism $\map h U U$ satisfying $h \cmp i \cmp e = j \cmp e$.
\end{pf}

Note that if $\id a$ is amalgamable in $\fK$ then this is indeed homogeneity (with respect to $a$). In particular, if $\fK$ has the amalgamation property then the weakly $\fK$-injective object is \define{homogeneous}{homogeneity}, that is, for every $\fL$-arrows $\map a i U$, $\map j a U$ with $a \in \ob{\fK}$ there exists an automorphism $\map h U U$ satisfying $h \cmp i = j$.
In general, the property of $U$ described in Theorem~\ref{THMigiuweg} can be called \emph{weak homogeneity}. We will elaborate this topic in the next section.

\section{Weak homogeneity}\label{SectSlabeHomogenita}

In the classical (model-theoretic) \fra\ theory, an important feature is that the \fra\ class can be reconstructed from its generic limit $U$, simply as the class of all finitely generated substructures (called the \emph{age} of $U$).
Actually, a countably generated model $U$ is the \fra\ limit of its age $\fK$ if and only if $U$ is homogeneous with respect to $\fK$, in the sense described above, where $\fK$ is treated as a category with embeddings.
That is why a \fra\ class is always assumed to be hereditary (i.e., closed under finitely generated substructures).
This cannot be formulated in category theory, however it becomes in some sense irrelevant, as we can always work in the category of \emph{all} finitely generated structures of a fixed language, or in a selected (usually full) subcategory.
On the other hand, we can consider subcategories of a fixed category $\fK$ and define the concept of being hereditary with respect to $\fK$.
By this way we can talk about objects that are weakly injective relative to a subcategory of $\fK$. We can also look at homogeneity and its weakening in a broader setting.

We continue using the framework from the previous section, namely, we assume that $\fK \subs \fL$ is a pair of categories satisfying (L0)--(L2).
Given a class of objects $\Ef \subs \ob{\fK}$, we can say that it is \define{hereditary}{hereditary class} in $\fK$ if for every $x \in \Ef$, for every $\fK$-arrow $\map f y x$ it holds that $y \in \Ef$.
Thus, the notion of being hereditary strongly depends on the category $\fK$ we are working with (the bigger category $\fL$ plays no role here).
Actually, it is more convenient, and within the philosophy of category theory, to define this concept for arbitrary subcategories (note that a class of objects may be viewed as a subcategory in which the arrows are precisely all the identities).
Namely, we say that a subcategory $\fS$ of $\fK$ is \define{hereditary}{hereditary subcategory} if for every compatible $\fK$-arrows $f,g$ the following equivalence holds:
$$f \cmp g \in \fS \iff f \in \fS.$$
Note that a hereditary subcategory $\fS$ is necessarily full. Indeed, if $\map f a b$ is such that $b \in \ob{\fS}$ then $\id b \in \fS$, therefore $f = \id b \cmp f \in \fS$.
It is straightforward to see that a family of objects $\Ef$ is hereditary if and only if the full subcategory $\fS$ with $\ob{\fS} = \Ef$ is hereditary, as a subcategory.
Conversely, if $\fS$ is a hereditary subcategory of $\fK$, then $\ob{\fS}$ is a hereditary class.

Natural examples of hereditary subcategories of $\fK$ are of the form
$$\fK_V := \setof{f \in \fK}{\fL(\cod(f),V) \nnempty},$$
where $V \in \ob{\fL}$.
One could call $\fK_V$ the \define{age}{age} of $V$ relative to $\fK$.
It is natural to ask when $\fK_V$ is a weak \fra\ category and when $V$ is its ``generic limit".
The answer is given below.

Fix $V \in \ob{\fL}$. We say that $V$ is \define{weakly homogeneous}{weakly homogeneous object} if for every $\fL$-arrow $\map f a V$ with $a \in \ob{\fK}$ there exist a $\fK$-arrow $\map e a b$ and an $\fL$-arrow $\map i b V$ such that $f = i \cmp e$ and for every $\fL$-arrow $\map j b V$ there is an automorphism $\map h V V$ satisfying $h \cmp f = j \cmp e$. This is shown in the following diagram in which, again, the triangle with vertices $b, V, V$ may not be commutative.
$$
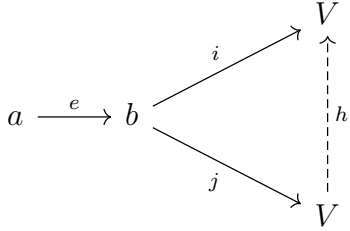
\begin{tikzcd}
& & & V  \ar[dd, dashrightarrow, "h"] \\
a \ar[rrru, bend left, "f"] \ar[r, "e"] & b \ar[rru, "i"] \ar[rrd, "j"'] & & \\
& & & V
\end{tikzcd}$$
Note that in this case, if $\map {j'}b V$ is another $\fL$-arrow, then there exists an automorphism $\map {h'} V V$ such that $h' \cmp f = j' \cmp e$. Thus, $k := h' \cmp h^{-1}$ is an automorphism of $V$ satisfying $k \cmp j \cmp e = j' \cmp e$.
This, by Theorem~\ref{THMigiuweg}, shows that the weakly $\fK$-injective object is weakly homogeneous. 

An $\fL$-object $V$ is \define{homogeneous}{homogeneous object} if the arrow $e$ in the definition above can always be identity. In other words, $V$ is homogeneous if for every $a\in \ob{\fK}$, for every $\fL$-arrows $\map i a V$, $\map j a V$ there is an automorphism $\map h VV$ satisfying $j = h \cmp i$. Homogeneity is often (especially by model-theorists) called \emph{ultra-homegeneity}.

The following result says that weakly homogeneous objects are weakly injective with respect to their age.

\begin{tw}\label{ThmSzescJedyn}
	Let $V \in \ob{\fL}$ and let $\fS := \fK_V$ be the age of $V$, as defined above. The following conditions are equivalent.
	\begin{enumerate}[itemsep=0pt]
		\item[{\rm(a)}] $V$ is weakly homogeneous.
		\item[{\rm(b)}] $\fS$ is a weak \fra\ category and $V$ is weakly $\fS$-injective.
	\end{enumerate}
\end{tw}

\begin{pf}
	(a)$\implies$(b)
	First, note that $V$ is $\fS$-cofinal.
	Fix an $\fL$-arrow $\map f a V$. Let $\map i b V$ and $\map e a b$ be as in the definition of weak homogeneity.
	Fix an arbitrary $\fS$-arrow $\map g b y$. There exists an $\fL$-arrow $\map k y V$. Apply the weak homogeneity to $j := k \cmp g$.
	By this way we obtain an automorphism $\map h V V$ satisfying $h \cmp k \cmp g \cmp e = f$.
	This shows that $V$ is weakly $\fS$-injective.

	Corollary~\ref{WNCorpietdwwwaa} says that $\fS$ is a weak \fra\ category (formally, one should replace $\fL$ by a suitable subcategory, so that (L1) will hold).
		
	(b)$\implies$(a) Trivial, by the comment after the definition of weak homogeneity.	
\end{pf}

Weak homogeneity was probably first studied by Pabion~\cite{Pab}, called \emph{prehomogeneity}, for multi-relations, i.e., structures with finitely many relations. This was later explored by Pouzet and Roux~\cite{PouRou}. We refer to~\cite{KraKub} for more details and bibliographic references.

We finish this section by exhibiting the (rather expected) relation between weak homogeneity and homogeneity. Recall that $\fK \subs \fL$ are as above, namely, conditions (L0)--(L2) are satisfied.

\begin{wn}
	Assume $V \in \ob \fL$ is weakly homogeneous.
	Then for every amalgamable $\fK$-arrow $\map e a b$, for every $\fL$-arrows $\map i a V$, $\map j b V$, there exists an automorphism $\map h V V$ such that $j \cmp e = h \cmp i \cmp e$.
	
	In particular, if $\fK$ has the amalgamation property, then $V$ is homogeneous.
\end{wn}

\begin{pf}
	Without loss of generality, we may assume $\fK = \fK_V$. Thus, $V$ is weakly $\fK$-injective and hence the statement follows directly from Theorem~\ref{THMigiuweg}. The second part is obvious, as the amalgamation property says that all identities are amalgamable.
\end{pf}

\section{The Banach-Mazur game}\label{SectBMGejmz}

In this section we explore connections between weakly injective objects and a natural infinite game which is a generalization of the classical Banach-Mazur game in topology.

We fix a category $\fK$. The \define{Banach-Mazur game}{Banach-Mazur game} played on $\fK$ is described as follows.
There are two players: \emph{Eve} and \emph{Odd}.
Eve starts by choosing $a_0 \in \ob{\fK}$.
Then Odd chooses $a_1 \in \ob{\fK}$ together with a $\fK$-arrow $\map {a_0^1}{a_0}{a_1}$.
More generally, after Odd's move finishing with an object $a_{2k-1}$, Eve chooses $a_{2k} \in \ob{\fK}$ together with a $\fK$-arrow $\map {a_{2k-1}^{2k}}{a_{2k-1}}{a_{2k}}$.
Next, Odd chooses $a_{2k+1} \in \ob{\fK}$ together with a $\fK$-arrow $\map {a_{2k}^{2k+1}}{a_{2k}}{a_{2k+1}}$.
Thus, the result of the play is a sequence
$$\map{\vec a}\omega \fK.$$
Of course, one needs to add the objective of the game, namely, a condition under which one of the players wins.
So, let us assume that $\fK$ is a subcategory of a bigger category $\fL$, so that \emph{some} sequences in $\fK$ have colimits in $\fL$.
For the moment, we do not need to assume neither of the conditions (L0)--(L2).
Now choose a family $\Wu \subs \ob \fL$.
We define the game $\BMG \fK \Wu$ with the rules described above, adding the statement that Odd \emph{wins the game} if and only if the colimit of the resulting sequence $\vec a$ is isomorphic to a member of $\Wu$.
So, Eve wins if either the sequence $\vec a$ has no colimit in $\fL$ or its colimit is isomorphic to none of the members of $\Wu$.

We are particularly interested in the case $\Wu = \sn W$ for some $W \in \ob{\fL}$, where the game $\BMG{\fK}{\Wu}$ will be denoted simply by $\BMG{\fK}{W}$.
Before we turn to it, we discuss some basic properties of our Banach-Mazur game.

Recall that a \define{strategy}{strategy} of Odd is a function $\Sigma$ assigning to each finite sequence $\map {\vec s} n \fK$ of odd length a $\fK$-arrow
$\map{\Sigma(\vec s)}{s_{n-1}}{s}$, called \emph{Odd's response} to $\vec s$.
We say that \emph{Odd plays according to} $\Sigma$ if the resulting sequence $\vec a$ satisfies $a_{n}^{n+1} = \Sigma(\vec a \rest n)$ for every odd $n \in \nat$.
Odd's strategy $\Sigma$ is \define{winning}{strategy!-- winning} in $\BMG{\fK}{\Wu}$ if $\lim \vec a$ is isomorphic to a member of $\Wu$ whenever Odd plays according to $\Sigma$, no matter how Eve plays.
These concepts are defined for Eve analogously.
A strategy $\Sigma$ of Eve is defined on sequences of even length, including the empty sequence, where $\Sigma(\emptyset)$ is simply a $\fK$-object $a_0$, the starting point of a play according to $\Sigma$.

\begin{tw}\label{ThmBMGweakdom}
	Let $\fK \subs \fL$ be two categories and let $\Wu \subs \ob{\fL}$.
	Let $\fS$ be a weakly dominating subcategory of $\fK$.
	Then Odd has a winning strategy in $\BMG{\fK}{\Wu}$ if and only if he has a winning strategy in $\BMG{\fS}{\Wu}$.
	The same applies to Eve.
\end{tw}

\begin{pf}
	Let $\Sigma$ be Odd's winning strategy in $\BMG{\fK}{\Wu}$.
	We describe his winning strategy in $\BMG{\fS}{\Wu}$.
	We denote the resulting sequence of a play in $\BMG{\fS}{\Wu}$ by $\vec s$.
	So, suppose Eve started with $s_0 \in \ob{\fS}$.
	Odd first chooses an $\fS$-arrow $\map {i_0}{s_0}{a_0}$ so that condition (W) of the definition of weak domination holds, namely, for every $\fK$-arrow $\map f{a_0}{x}$ there is a 
	$\fK$-arrow $\map g x t$ such that $g \cmp f \cmp i_0 \in \fS$.
	Let $a_0^1 = \Sigma(a_0)$, so $\map {a_0^1}{a_0}{a_1}$ with $a_1 \in \ob{\fK}$.
	Using (W), Odd finds a $\fK$-arrow $\map {j_0}{a_1}{s_1}$ and he responds with $s_0^1 := j_0 \cmp a_0^1 \cmp i_0$.
	In general, the strategy is described in the following commutative diagram.
	$$\begin{tikzcd}
	s_0 \ar[r] \ar[d, "i_0"'] & s_1 \ar[r] & \cdots \ar[r] & s_{2n} \ar[d, "i_n"'] \ar[rr] & & s_{2n+1} \ar[r] & \cdots \\
	a_0 \ar[r, "a_0^1"'] & a_1 \ar[u, "j_0"'] & & a_{2n} \ar[rr, "a_{2n}^{2n+1}"'] & & a_{2n+1} \ar[u, "j_n"']
	\end{tikzcd}$$
	Namely, when Eve finishes with $s_{2n}$, Odd first chooses a suitable $\fS$-arrow $\map{i_n}{s_{2n}}{a_{2n}}$ realizing the weak domination. Next, he uses $\Sigma$ to find a $\fK$-arrow $\map{a_{2n}^{2n+1}}{a_{2n}}{a_{2n+1}}$.
	Specifically, $a_{2n}^{2n+1}$ is Odd's response to the sequence
	$a_0 \to a_1 \to \dots \to a_{2n}$
	in which the arrows are suitable compositions of those from the diagram above.
	Odd responds with $s_{2n}^{2n+1} := j_n \cmp f_n \cmp i_n$, where $j_n$ comes from the weak domination of $\fS$ (condition (W)).
	This is a winning strategy, because the resulting sequence $\vec s$ is isomorphic to the sequence $\vec a$, where
	$$a_{2k+1}^{2k+2} = i_{k+1} \cmp s_{2k+1}^{2k+2} \cmp j_k$$
	for every $k \in \nat$; this sequence is the result of a play of $\BMG{\fK}{\Wu}$ in which Odd was using strategy $\Sigma$.
	
	Now suppose Odd has a winning strategy $\Sigma$ in $\BMG{\fS}{\Wu}$.
	Playing the game $\BMG{\fK}{\Wu}$, assume Eve started with $a_0 \in \ob{\fK}$.
	Odd first uses (C) to find an arrow $\map{i_0}{a_0}{s_0}$ with $s_0 \in \ob{\fS}$.
	Next, he takes the arrow $\map {s_0^1}{s_0}{s_1}$ according to $\Sigma$. Specifically, $s_0^1 = \Sigma(s_0)$. He responds with $a_0^1 := j_0 \cmp s_0^1 \cmp i_0$, where $\map{j_0}{s_1}{a_1}$ is an $\fS$ from condition (W), namely, for every $\fK$-arrow $\map f{a_1}x$ there is a $\fK$-arrow $\map g x s$ satisfying $g \cmp f \cmp j_0 \in \fS$.
	In general, the strategy described in the following commutative diagram.
	$$\begin{tikzcd}
	a_0 \ar[d, "i_0"'] \ar[r] & a_1 \ar[r] & \cdots \ar[r] & a_{2n} \ar[d, "i_n"'] \ar[rr] & & a_{2n+1} \ar[r] & \cdots \\
	s_0 \ar[r, "s_0^1"'] & s_1 \ar[u, "j_0"'] & & s_{2n} \ar[rr, "s_{2n}^{2n+1}"'] & & s_{2n+1} \ar[u, "j_n"']
	\end{tikzcd}$$
	Here, $i_n$ comes from condition (W), namely, $i_n \cmp a_{2n-1}^{2n} \cmp j_{n-1} \in \fS$.
	Furthermore, $s_{2n}^{2n+1} = \Sigma(\vec v)$, where $\vec v$ is the sequence
	$s_0 \to s_1 \to s_2 \to \dots \to s_{2n}$
	obtained from the diagram above (note that all its arrows are in $\fS$).
	Finally, $j_n$ is such that the assertion of (W) holds, that is, for every $\fK$-arrow $\map{f}{a_{2n+1}}{x}$ there is a $\fK$-arrow $\map g x t$ such that $g \cmp f \cmp j_n \in \fS$.
	Odd's response is $a_{2n}^{2n+1} := j_n \cmp s_{2n}^{2n+1} \cmp i_n$.
	This strategy is winning in $\BMG{\fK}{\Wu}$, because the resulting sequence $\vec a$ is isomorphic to the sequence $\vec s$ in which
	$$s_{2k+1}^{2k+2} = i_{k+1} \cmp a_{2k+1}^{2k+2} \cmp j_k \in \fS$$
	for every $k \in \nat$.
	The sequence $\vec s$ results from a play of $\BMG{\fS}{\Wu}$ in which Odd was using his winning strategy $\Sigma$.
	
	The case of Eve's winning strategies is almost the same, as the rules are identical for both players, except for Eve's first move.
\end{pf}

%The next result can be viewed as a category-theoretic version of the Baire Category Theorem.

\begin{tw}\label{ThmBaireCatThC}
	Assume $\sett{\Wu_n}{\ntr}$ is such that each $\Wu_n \subs \ob \fL$ is closed under isomorphisms and Odd has a winning strategy in $\BMG \fK {\Wu_n}$ for each $\ntr$.
	Then Odd has a winning strategy in
	$$\BMG \fK {\bigcap_{\ntr}\Wu_n}.$$
	In particular, $\bigcap_{\ntr}\Wu_n \nnempty$.
\end{tw}

\begin{pf}
	Let $\Sigma_n$ denote Odd's winning strategy in $\BMG{\fK}{\Wu_n}$.
	Let $\ciag{I}$ be a partition of all even natural numbers into infinite sets.
	Let $J_n = I_n \cup \setof{i+1}{i \in I_n}$.
	Given a finite sequence $\vec s$ whose length $n$ is odd, let $k$ be such that $n-1 \in I_k$ and define
	$$\Sigma(\vec s) = \Sigma_k(\vec s \rest (J_k \cap n)).$$
	We claim that $\Sigma$ is a winning strategy of Odd in the game $\BMG{\fK}{\bigcap_\ntr \Wu_n}$.
	
	Indeed, suppose $\vec a$ is the result of a play in which Odd has been using strategy $\Sigma$.
	Then $\vec a \rest J_k$ is a sequence resulting from another play in which Odd was using strategy $\Sigma_k$.
	Thus $\lim \vec a = \lim(\vec a \rest J_k) \in \Wu_k$. Hence $\lim \vec a \in \bigcap_\ntr \Wu_n$.
\end{pf}

We now switch to the case where $\Wu$ is the isomorphism class of a single object.
As the reader may guess, weakly injective objects play a significant role here.
In the next result we do not assume (L0)--(L2).

\begin{tw}\label{ThmGnrikOdws}
	Let $\fK \subs \fL$ and assume that $\vec u$ is a weak \fra\ sequence in $\fK$ with $U = \lim \vec u$ in $\fL$.
	Then Odd has a winning strategy in $\BMG{\fK}{U}$.
\end{tw}

\begin{pf}
	We may assume that the sequence $\vec u$ is normalized.
	Odd's strategy is as follows.
	Suppose $a_0 \in \ob{\fK}$ is Eve's first move. Using (G1), Odd finds $k \in \nat$ together with a $\fK$-arrow $\map{f_0}{a_0}{u_{k}}$. His response is $a_0^1 := u_{k}^{k+1} \cmp f_0$. In particular, $a_1 = u_{k+1}$.
	
	In general, suppose $a_{2n-1}^{2n}$ was the $n$th move of Eve. Assume inductively that $a_{2n-1} = u_{\ell+1}$ and $a_{2n-2}^{2n-1} = u_{\ell}^{\ell+1} \cmp f_{n-1}$ for some $\fK$-arrow $f_{n-1}$.
	Using (G2), Odd finds $m > \ell+1$ together with a $\fK$-arrow $\map{f_n}{a_{2n}}{u_{m}}$ satisfying
	$$u_{\ell}^{m} = f_n \cmp a_{2n-1}^{2n} \cmp u_{\ell}^{\ell+1}.$$
	Odd's response is
	$a_{2n}^{2n+1} := u_{m}^{m+1} \cmp f_n$.
	In particular, $a_{2n+1} = u_{m+1}$. The strategy is shown in the following diagram.
	$$\begin{tikzcd}
	\cdots \ar[r] & u_{\ell} \ar[r] & u_{\ell+1} \ar[rd, "a_{2n-1}^{2n}"'] & & u_{m}  \ar[r] & u_{m+1} \ar[r] & \cdots \\
	& & & a_{2n} \ar[ru, "f_n"']
 	\end{tikzcd}$$
	It is clear that the resulting sequence $\vec a$ is isomorphic to $\vec u$, therefore $\lim \vec a = U$.
\end{pf}

The proof above is somewhat similar to that of Theorem~\ref{ThmBMGweakdom}.
In fact, if the sequence $\vec u$ is one-to-one (that is, $u_n \ne u_m$ for $n \ne m$) then one can use Theorem~\ref{ThmBMGweakdom} to play the game in the image of $\wek u$, where Odd's winning strategy is obvious.

Our goal is reversing Theorem~\ref{ThmGnrikOdws}, extending the results of Krawczyk and the author~\cite{KraKub}.
We start with a technical lemma.
Recall that a category $\fC$ is \emph{locally countable} if $\fC(x,y)$ is a countable set for every $\fC$-objects $x,y$.

\begin{lm}\label{LMxEwa}
	Assume $\fK \subs \fL$ are two categories satisfying (L0)--(L2), $\fK$ is locally countable, $V \in \ob{\fL}$, and suppose $\map e a V$ is an $\fL$-arrow with $a \in \ob{\fK}$ satisfying the following condition.
	\begin{enumerate}
		\item[{\rm($\times$)}] For every $\fK$-arrow $\map f a b$ there exists a $\fK$-arrow $\map{f'}{b}{b'}$ such that for every $\fL$-arrow $\map i{b'}V$ it holds that $e \ne i \cmp f' \cmp f$.
	\end{enumerate}
	Then Eve has a winning strategy in $\BMG{\fK}{V}$.
\end{lm}

\begin{pf}
	Eve's strategy is as follows.
	She starts with $a_0 := a$.
	At step $n>0$, Eve chooses a $\fK$-arrow $\map{f_n}a{a_{2n-1}}$
	and responds with $a_{2n-1}^{2n} := f'$, where $f'$ comes from condition ($\times$) applied to $f := f_n$.
	Thus
	\begin{equation}
		(\forall\; i \in \fL(a_{2n}, V)) \;\; e \ne i \cmp a_{2n-1}^{2n} \cmp f_n.
		\tag{$\dagger$}\label{Eqsztyletjeden}
	\end{equation}
	Of course, this strategy depends on the choice of the sequence $\sett{f_n}{n>0}$. We show that a suitable choice makes Eve's strategy winning.
	Namely, she needs to take care of all $\fK$-arrows from $a$ into the sequence $\vec a$.
	More precisely, the following condition should be satisfied.
	\begin{equation}
		(\forall\; k>0) (\forall\; g \in \fK(a, a_k)) (\exists\; n>k) \;\; f_n = a_k^{2n-1} \cmp g.
		\tag{$\ddagger$}\label{EqDdaggerBookkeepin}
	\end{equation}
	In order to achieve (\ref{EqDdaggerBookkeepin}), we use the fact that $\fK$ is locally countable.
	Specifically, for each $k>0$, for each $g \in \fK(a, a_k)$ we inductively choose an integer $\phi(k,g) > k$ in such a way that $\phi(k',g') \ne \phi(k,g)$ whenever $\pair k g \ne \pair{k'}{g'}$. This is possible, because for a fixed $k$ there are only countably many possibilities for $g$ (we may first partition $\omega$ into infinite sets $B_k$ and make sure that $\phi(k,g) \in B_k$ for every $g$).
	We set $f_n := a_k^{2n-1} \cmp g$ whenever $n = \phi(k,g)$.
	
	Now let $A = \lim \vec a \in \ob{\fL}$ and suppose that $\map h V A$ is an isomorphism in $\fL$.
	Using (L2), we find a $\fK$-arrow $\map g a {a_k}$ such that $h \cmp e = a_k^\infty \cmp g$, where $a_k^\infty$ is part of the colimiting co-cone.
	By (\ref{EqDdaggerBookkeepin}), there is $n > k$ such that $f_n = a_k^{2n-1} \cmp g$.
	Consider $i := h^{-1} \cmp a_{2n}^\infty$. We have
	$$i \cmp a_{2n-1}^{2n} \cmp f_n = h^{-1} \cmp a_{2n}^\infty \cmp a_{2n-1}^{2n} \cmp a_k^{2n-1} \cmp g = h^{-1} \cmp a_k^\infty \cmp g = h^{-1} \cmp h \cmp e = e,$$
	contradicting (\ref{Eqsztyletjeden}).
	This shows that Eve wins while using the strategy above.
\end{pf}

We are ready to prove the main result of this section.

\begin{tw}\label{ThmSedmPet}
	Assume $\fK \subs \fL$ satisfy (L0)--(L2) and $\fK$ is locally countable.
	Given an $\fL$-object $V$, the following properties are equivalent.
	\begin{enumerate}[itemsep=0pt]
		\item[{\rm(a)}] $V$ is weakly $\fK$-injective (in particular, $\fK$ is a weak \fra\ category).
		\item[{\rm(b)}] Odd has a winning strategy in $\BMG{\fK}{V}$.
		\item[{\rm(c)}] Eve does not have a winning strategy in $\BMG{\fK}{V}$.
	\end{enumerate}
\end{tw}

\begin{pf}
	(a)$\implies$(b)
	By (L1), $V = \lim \vec v$ for a sequence $\vec v$ in $\fK$. By Theorem~\ref{ThmMgenerikgnrc}, this sequence is weak \fra\ in $\fK$.
	Thus (b) follows from Theorem~\ref{ThmGnrikOdws}.
	
	(b)$\implies$(c) Obvious.
	
	(c)$\implies$(a)
	First, note that $V$ satisfies (U), since if $x \in \ob{\fK}$ is such that $\fL(x, V) = \emptyset$ then Eve would have an obvious winning strategy, starting the game with $x$.
	Thus, supposing $V$ is not weakly $\fK$-injective, we deduce that it is not weakly $\fK$-injective.
	Hence, there exists $\map e a V$ with $a \in \ob{\fK}$ such that for every $\fK$-arrow $\map f a b$ there is a $\fK$-arrow $\map {f'} b y$ such that no $\fL$-arrow $\map j y V$ satisfies $j \cmp f' \cmp f = e$.
	This is precisely condition ($\times$) of Lemma~\ref{LMxEwa}, contradicting (c).
\end{pf}

Note that the result above says, in particular, that the Banach-Mazur game played on a locally countable category is determined, as long as the goal is a single isomorphic type.
This is not true when the goal is an arbitrary set of objects, as the next example shows.

\begin{ex}
	Let $X \ne \emptyset$ be a compact Hausdorff topological space, $\fK$ the family of all nonempty open subsets of $X$.
	Define $U \rloe V$ if either $U=V$ or $\cl V \subs U$, where $\cl V$ denotes the closure of $V$.
	Then $\pair \fK \rloe$ is a poset, therefore it is a category.
	Let $\fL$ be the family of all nonempty closed $G_\delta$ subsets of $X$, endowed with the same ordering.
	Then $\fK \subs \fL$ satisfies (L0)--(L2). Clearly, $\fK$ is locally countable, being a poset category. It is also clear that $\fK$ fails the WAP, unless $X$ is a singleton.
	Now the Banach-Mazur game played on $\fK$ is practically the same as the original topological Banach-Mazur game, the only difference is that we force a stronger containment relation, in order to achieve (L2).
	If, additionally, $X$ is a metric space then each of the players can play so that the intersection of the resulting sequence is a single point.
	Thus, the game $\BM{\fK}$ can be parameterized by a subset $Y$ of $X$, meaning that Odd wins in $\BMG{\fK}{Y}$ if the intersection of the resulting chain is an element of $Y$.
	Now, if $Y$ is a Bernstein set then none of the players has a winning strategy in $\BMG{\fK}{Y}$.
	Recall that a set $Y$ is \define{Bernstein}{Bernstein set} if for every perfect set $P \subs X$ it holds that $P\cap Y \nnempty \ne P \setminus Y$.
	A well known characterization due to Oxtoby~\cite{Oxt} says that Odd has a winning strategy in $\BMG{\fK}{Y}$ if and only if $Y$ contains a set $G$ that is both dense and $G_\delta$ in $X$.
	On the other hand, Eve has a winning strategy if and only if there is a nonempty open set $U \subs X$ such that $U \setminus Y$ contains a $G_\delta$ set whose closure contains $U$. Every nonempty $G_\delta$ set contains a perfect set, therefore if $Y$ is Bernstein then none of the players can have a winning strategy.
\end{ex}

The example above shows also that the abstract Banach-Mazur game indeed generalizes the classical topological one, invented by Mazur around 90 years ago.

Motivated by the results above, we now introduce the concept of a generic object. In the context of model theory, it already appeared in \cite{KubBMG}.

\begin{df}[Generic object]
	Let $\fK \subs \fL$ be as above, satisfying (L0)--(L2).
	An object $V \in \ob{\fL}$ will be called \define{$\fK$-generic}{generic object} if Odd has a winning strategy in the game $\BMG{\fK}{V}$.
\end{df}

Theorem~\ref{ThmSedmPet} tells us that this concept coincides with weak injectivity, as long as the category is locally countable.
On the other hand, the above definition covers natural categories that are not locally countable, e.g., the category of all finite-dimensional normed spaces with linear isometric embeddings, which apparently has a generic object in the category of all separable Banach spaces (see \cite{KGur}). In any case, the definition of a generic object is more general and perhaps more natural than that of a weakly injective object.

Let us admit that the adjective \emph{generic} has already been used by Kueker and Laskowski~\cite{KueLas} referring to model-theoretic \fra\ limits where the notion of embedding is specialized---from the point of view of category theory this is just selecting a wide subcategory of the category of all embeddings.

\subsection{Ubiquity of generic objects}

In model theory, the question which isomorphic types of objects can be represented as a residual set in a suitable complete metric space was addressed by Cameron~\cite{Cam} and later explored by Pouzet and Roux~\cite{PouRou} who actually proved that this happens if and only if the object is weakly homogeneous (pre-homogeneous, in their terminology).
The results concerning the abstract Banach-Mazur game, combined with Oxtoby's characterization of winning strategies lead to a direct proof of Pouzet and Roux' result saying that weakly homogeneous objects are ubiquitous in the sense of Baire category\footnote{Here, the word \emph{category} comes from the Baire Category Theorem.}.

In order to present this result, we need to define a suitable complete metric space.
Namely, let $\fK \subs \fL$ be as before, satisfying (L0)--(L2) and assume that $\fK$ is small. In most cases, we may replace $\fK$ by a small (often countable) subcategory, simply by localizing the objects in a big enough set.

Now, let $T_\fK = \fK^{<\omega}$ be the set of all finite sequences in $\fK$ endowed with the ``end-extension" ordering. Then $X$ is a tree and each of its branches corresponds to an object of $\fL$ (namely, its colimit).
Let $\bd T_\fK$ be the space of all branches through $T_\fK$.
Then $X$ is a complete ultrametric space, when endowed with the metric $\rho(x,y) = 1/n$, where $n$ is minimal with the property $x^{n}_{n_1} \ne y^n_{n-1}$
(it is tempting to replace this by $x_n \ne y_n$, however the category $\fK$ may consist of a single object and then only the arrows distinguish the branches of $T_\fK$).
In this setting, our Banach-Mazur game is equivalent to the classical topological Banach-Mazur game played with basic open sets (which are actually ultra-metric balls). Thus, Oxtoby's characterization applies. Summarizing:

\begin{tw}
	Let $T_\fK$ and $\bd T_\fK$ be as above and $\fK$ is locally countable.
	If $\fK$ has a generic object $U$ in $\fL$ then the set of all $v \in \bd T_\fK$ with $\lim v \iso U$ is dense $G_\delta$ in $\bd T_\fK$.
\end{tw}

\begin{pf}
	The fact that this set is dense follows directly from Oxtoby's result~\cite{Oxt}. In order to show that it is $G_\delta$, it suffices to observe that $\fK$ is weakly dominated by a countable subcategory and that the definition of a weak \fra\ sequence requires countably many parameters, all of them defining open subsets of $\bd T_\fK$.	
\end{pf}

\subsection{Cofinality vs. WAP}

We now extend the results from~\cite{KraKub} concerning the weak amalgamation property.
Namely, let us fix a pair of categories $\fK \subs \fL$ and for the moment the only assumption we make is that $\fK$ be directed.
We define the \emph{cofinality number} as follows:
$$\cf(\fL) := \min \setof{|\Yu|}{\Yu \text{ is cofinal in } \fL},$$
where a family of $\fL$-objects $\Yu$ is \emph{cofinal} if for every $\fL$-object $X$ there exists an $\fL$-arrow from $X$ to some $U \in \Yu$.
In model theory, cofinality is called \emph{universality}, however in category theory there would be a conflict with the notion of a universal object, which has a different meaning.

In any case, if $\fK$ is a \fra\ category and (L0)--(L2) hold, then $\cf(\fL) = 1$, which is witnessed by the \fra\ limit $U \in \ob{\fL}$.
On the other hand, there are examples of weak \fra\ categories of finite graphs $\fK$, where $\fL$ is the corresponding category of countable graphs, such that $\cf(\fL) = \cont$, the cardinality of the set of all reals $\Err$ (the continuum), see~\cite{KraKub}.

On the other hand, failure of the weak amalgamation property implies that the cofinality number is large:

\begin{tw}\label{THMbweigwieug}
	Assume $\fK \subs \fL$ are categories satisfying (L0)--(L2), $\fK$ is directed and locally countable. If $\cf(\fL) < \cont$ then $\fK$ has the weak amalgamation property.
\end{tw}

Before proving the theorem above, we formulate a result involving the Banach-Mazur game, that could be of independent interest.
Given $W \in \ob{\fL}$, denote by $\below{W}$ the class of all $X \in \ob{\fL}$ for which $\fL(X, W) \nnempty$.
Similarly, if $\Yu$ is a class of $\fL$-objects, we denote
$$\below{\Yu} = \bigcup_{W \in \Yu} \below{W}.$$
Thus, the second player wins the game $\BMG{\fK}{\below{\Yu}}$
if and only if the colimit of the resulting sequence admits at least one $\fL$-arrow into some $W \in \Yu$.
Theorem~\ref{THMbweigwieug} is an immediate consequence of the following fact.

\begin{prop}
	Assume $\fK \subs \fL$ are categories satisfying (L0)--(L2), $\fK$ is directed and locally countable. If Odd has a winning strategy in $\BMG{\fK}{\below{\Yu}}$ for some family $\Yu \subs \ob{\fL}$ of cardinality $< \cont$, then $\fK$ has the weak amalgamation property.
\end{prop}

\begin{pf}
	Let $\cantree$ denote the tree of all finite zero-one sequences.
	Given $s \in \cantree$ and $i \in 2$, we denote by $s \concat i$ the sequence obtained from $s$ by adding $i$ at the end. Given a sequence $s \in \cantree$ of positive length (the only sequence of length zero is $\emptyset$), we denote by $s^-$ the sequence obtained from $s$ by removing its last element. So if the last element was $i \in 2$, then $s = (s^-)\concat i$.
	
	Suppose $z$ is a witness for the failure of WAP.
	Using induction, we build a family of $\fK$-objects $\sett{z_s}{s \in \cantree}$ together with $\fK$-arrows $\map{e_s}{z_{s^-}}{z_s}$, such that the following conditions are satisfied.
	\begin{enumerate}[itemsep=0pt]
		\item[(1)] $z_\emptyset = z$.
		\item[(2)] The arrows $e_{s{\concat0}}$, $e_{s{\concat1}}$ witness the fact that the arrow 
		$$\map{e_{s}\cmp \dots \cmp e_{s \rest 2}\cmp e_{s \rest 1}}{z_\emptyset}{z_s}$$
		is not amalgamable.
		\item[(3)] For each $\sig \in \Cantor$, the sequence
		$$\begin{tikzcd}
			z_\emptyset \ar[r] & z_{\sig \rest 1} \ar[r] & \cdots \ar[r] & z_{\sig \rest n} \ar[r] & \cdots
		\end{tikzcd}$$
		comes from a play according to Odd's winning strategy in $\BMG{\fK}{\below{\Yu}}$.
	\end{enumerate}
	It is clear how to achieve (2), knowing that $z_\emptyset$ witnesses the failure of the WAP. In order to achieve (3), we need to add the following clause:
	\begin{enumerate}
		\item[(4)] For each $s \in \cantree$, for each $i \in 2$ there are $z_{s\concat i}'$ and $\fK$-arrows $\map{e_{s\concat i}'}{z_s}{z_{s\concat i}'}$, $\map{e_{s\concat i}''}{z_{s\concat i}'}{z_{s\concat i}}$ such that $e_{s\concat i} = e_{s\concat i}'' \cmp e_{s\concat i}'$ and $e_{s\concat i}''$ is Odd's answer, according to his winning strategy, to the sequence
		$$\begin{tikzcd}
			z_\emptyset \ar[r, "e_{s \rest 1}'"] & z_{s \rest 1}' \ar[r, "e_{s \rest 1}''"] & z_{s \rest 1} \ar[r, "e'_{s \rest2}"] & z'_{s \rest2} \ar[r, "e''_{s \rest2}"] & z_{s \rest2} \ar[r] & \cdots \ar[r] & z_{s} \ar[r, "e_{s \concat i}'"] & z_{s \concat i}'
		\end{tikzcd}$$
		where the first move of Eve is actually $z_{s \rest 1}'$ instead of $z_\emptyset$.
	\end{enumerate}
	To be more precise, the inductive step at $s \in \cantree$ runs as follows. We first choose $e_{s \concat0}'$, $e_{s \concat1}'$ as in (2) and then we compose them with $e_{s \concat0}''$, $e_{s \concat1}''$, using Odd's strategy, as described in (4).
	
	Finally, for each $\sig \in \Cantor$ we have a branch of our tree, namely, a sequence $\map{\vec{z_\sig}}{\nat}{\fK}$ such that $\vec{z_\sig}(n) = z_{\sig \rest n}$ and the bonding arrow from $\vec{z_\sig}(n-1)$ to $\vec{z_\sig}(n)$ is $e_{\sig \rest n}$.
	All these sequences are results of instances of playing the game $\BMG{\fK}{\below{\Yu}}$, where Odd was using his winning strategy.
	Denote by $e_\sig$ the colimiting arrow from $z = z_\emptyset$ to the colimit of $\vec{z_\sig}$.
	Choose $U_\sig \in \Yu$ so that there is $ i_\sig \in \fL(\lim \vec{z_\sig}, U_\sig)$. Since $|\Yu| < \cont$, there are $W \in \Yu$ and an uncountable set $S \subs \Cantor$ such that $U_\sig = W$ for every $\sig \in S$.
	
	Using (L1), (L2) and the fact that $\fK$ is locally countable, we can find two (in fact, uncountably many) different $\sig, \tau \in S$ such that
	$$i_\sig \cmp e_\sig = i_\tau \cmp e_\tau.$$
	Let $s = \sig \meet \tau \in \cantree$ be the maximal common part of $\sig$ and $\tau$.
	We may assume $s \concat0 \subs \sig$ and $s \concat1 \subs \tau$.
	Let $e_{\sig \setminus s}$ denote the colimiting arrow from $z_{s \concat0}$ to $\lim \vec{z_\sig}$ and let $e_{\tau \setminus s}$ denote the colimiting arrow from $z_{s \concat1}$ to $\lim \vec{z_\tau}$.
	Using (L1), we find a sequence $\vec w$ in $\fK$ whose colimit is $W$.
	Using (L2), we see that both  
	$i_\sig \cmp e_{\sig \setminus s}$, $i_\tau \cmp e_{\tau \setminus s}$ factor through a fixed $w_n$, that is,
	$$i_\sig \cmp e_{\sig \setminus s} = w_n^\infty \cmp f_0 \oraz 
	i_\tau \cmp e_{\tau \setminus s} = w_n^\infty \cmp f_1$$
	for some $\map{f_0}{z_{s \concat0}}{w_n}$, $\map{f_1}{z_{s \concat1}}{w_n}$. The situation is described in the diagram below, where $k = e_{s}\cmp \dots \cmp e_{s \rest 2}\cmp e_{s \rest 1}$.
	$$\begin{tikzcd}
		&&& {z_{s \concat1}} && {} & {\lim \vec{z_\tau}} \\
		{z_{\emptyset}} && {z_s} &&& 	{w_n} &&& W \\
		&&& {z_{s \concat0}} && {} & 	{\lim \vec{z_\sigma}}
		\arrow["k"{description}, from=2-1, to=2-3]
		\arrow["{w_n^\infty}"{description}, from=2-6, to=2-9]
		\arrow["e_{s \concat1}", from=2-3, to=1-4]
		\arrow["e_{s \concat0}"', from=2-3, to=3-4]
		\arrow["{e_\tau}", bend left, from=2-1, to=1-7]
		\arrow["{e_\sigma}"', bend right, from=2-1, to=3-7]
		\arrow["{i_\tau}", from=1-7, 	to=2-9]
		\arrow["{i_\sigma}"', from=3-7, to=2-9]
		\arrow["{e_{\tau \setminus s}}"{description}, from=1-4, to=1-7]
		\arrow["{e_{\sigma \setminus 	s}}"{description}, from=3-4, to=3-7]
		\arrow["{f_1}", from=1-4, to=2-6]
		\arrow["{f_0}"', from=3-4, to=2-6]
	\end{tikzcd}$$
	Finally, using the fact that $w_n^\infty$ is a monic (condition (L0)), we obtain that $f_0 \cmp e_{s \concat0} \cmp k = f_1 \cmp e_{s \concat1} \cmp k$, contradicting the fact that $k$ is not amalgamable (condition (2)). This completes the proof.
\end{pf}

The result above says that we can deduce the weak amalgamation property from the small cofinality number.
On the other hand, one can find in~\cite{KraKub} the following relevant example, showing that the converse is false: $\fK$ is the class of all finite graphs with vertex degree $\loe3$ and $\fL$ is the category of all countable graphs of vertex degree $\loe3$. Then $\fK \subs \fL$ satisfy conditions (L0)--(L2), $\fK$ is directed, essentially countable, and has the cofinal amalgamation property, while $\cf(\fL) =\cont$. Another example in~\cite{KraKub} is a hereditary class of finite graphs without the WAP.
We close this section with two simple examples of countable directed categories failing the WAP. Both of them share the same idea, even though they lie on the opposite sides of category theory (one is a monoid and the other one comes from a poset).

\begin{ex}\label{EXnefEFNFq}
	Let $\fK$ be the free monoid over a countable infinite set of letters, say, $\nat$. So $\fK$ consists of all words of the form
	$$n_0 n_1 \dots n_{k-1},$$
	where $\sett{n_i}{i<k} \subs \nat$ and the empty word is the identity.
	Clearly, $\fK$ treated as a category is countable and directed (it has just one object). On the other hand, it obviously fails the weak amalgamation property.
	Note that any category $\fL \sups \fK$ satisfying (L0)--(L2) must have cofinality $\goe \cont$. A natural choice is the category of all countable words in the alphabet $\nat$. An arrow from a countable word $x$ to another countable word $y$ is uniquely determined by cutting off initial parts of $x$ and $y$ in the sense that $x$ and $y$ are equal modulo a finite initial part.
\end{ex}

\begin{ex}\label{ExoHIEF}
	Let $\fK$ be the free category over the poset $\pair{\nat}{\loe}$ treated as a directed graph.
	So the objects of $\fK$ are natural numbers and the arrows are paths that could be encoded as finite increasing sequences of natural numbers.
	More precisely, an arrow from $k$ to $\ell$ is any sequence of the form
	$$k= m_0 < m_1 < \dots < m_{n-1} \loe m_n = \ell.$$
	We allow $m_{n-1} = m_n$ in order to get the identities, namely, $\pair \ell \ell$ is the identity of $\ell$. It is rather obvious that $\fK$ fails the weak amalgamation property. On the other hand, $\fK$ is countable and directed. Note also that all $\fK$-arrows are monic .
	
	What is a natural category $\fL \sups \fK$ satisfying (L0)--(L2)?
	A sequence in $\fK$ can be identified with a finite or infinite increasing sequence of natural numbers, namely, a unique increasing enumeration of a subset of $\nat$.
 	It is easy to check that an arrow from a sequence $\vec x$ to a sequence $\vec y$ exists if and only if the corresponding sets $X$, $Y$ are equal, modulo a finite set.
 	We conclude that a natural choice of $\fL$ is the category of all subsets of $\nat$ with suitable arrows (namely, inclusions modulo a finite set).
 	Note that $\fL$ has cofinality $\cont$, due to Theorem~\ref{THMbweigwieug}.
\end{ex}

\section{Applications}\label{SectAppsyJ}

The theory presented above definitely calls for illustrative examples. Some of them have already been described above. Many relevant examples can be found in the existing literature (e.g. \cite{Pab}, \cite{Ivanov}, \cite{KechRos}, \cite{KraKub}). Below we collect some groups of examples, focusing on the weak amalgamation property, which is in fact the main theme of this note.

\subsection{Monoids}

Category theory has two extremes: monoids and posets, living on the opposite sides.
Thus, before going into concrete categories of models, it is natural to discuss the theory of weak \fra\ categories in those extreme settings.
It turns out that posets do not contribute much, as all diagrams in a poset category are commutative, therefore WAP is equivalent to AP.
A poset is a \fra\ category if and only if it is directed and has a countable cofinality.
On the other hand, the free category over a quite simple poset already leads to something nontrivial, as shown in Example~\ref{ExoHIEF} above. Furthermore, some countable monoids fail the WAP, see Example~\ref{EXnefEFNFq} above.

Let $\fM = \triple{M}{\cmp}{1}$ be a monoid, treated as a category with the unique object $M$.
Note that if $\fM$ is commutative then it obviously has the amalgamation property.
A countable monoid without the WAP has already been described in Example~\ref{EXnefEFNFq}.

The weak amalgamation property in $\fM$ can be rephrased as follows: There exists $e \in M$ such that for every $x,y \in M$ there are $x',y' \in M$ satisfying
$$x' \cmp x \cmp e = y' \cmp y \cmp e.$$
We are particularly interested in left-cancellative monoids (as condition (L0) suggests). It is well known (and easy to show) that every left-cancellative monoid is isomorphic to a monoid of one-to-one transformations of a fixed set $S$. We may assume that $S$ has an extra structure and the elements of the monoid are embeddings or just one-to-one maps preserving the structure.
Now it is easy to see examples of non-commutative monoids with the amalgamation property.
Perhaps the simplest one is the monoid of all one-to-one self-maps of a fixed infinite set $S$.
A slightly more sophisticated example is given below.

\begin{ex}
	Let $R$ be a ring, $n \goe 1$, and let $\fM$ be the monoid of all $(n \times n)$-matrices with coefficients in $R$ and with non-zero determinant. The operation is multiplication.
	If $n>1$ then $\fM$ is not commutative. Assume $R$ is an integral domain and let $K$ be its field of fractions.
	Then $\fM$ has the amalgamation property.
	
	Indeed, given $A, B \in \fM$, we have matrices $A^{-1}, B^{-1}$ with coefficients in $K$. Let $g \in R$ be the product of all denominators of the fractions appearing in $A^{-1}$ and $B^{-1}$. Then $gA^{-1}$ and $gB^{-1}$ are matrices with coefficients in $R$, therefore they are elements of $\fM$.
	Finally, $(gA^{-1})A = gI = (gB^{-1})B$, where $I$ is the identity matrix.
	
	Actually, even $R := \Zee$ with $n := 1$ leads to a nontrivial example. Namely, this corresponds to the monoid of all self-embeddings of the group $\pair{\Zee}{+}$. Thus, in this case the natural category $\fL \sups \fM$ satisfying conditions (L0)--(L2) is the category of all countable abelian groups whose all nontrivial finitely generated subgroups are isomorphic to $\pair{\Zee}{+}$. It is easy to see that the \fra\ limit of $\fM$ in $\fL$ is the group $\pair{\Qyu}{+}$.
	It is also easy to check that if $n \goe 1$ then $\pair{\Qyu^n}{+}$ is the \fra\ limit, where now the matrices encode self-embeddings of $\pair{\Zee^n}{+}$.
\end{ex}	

The following example, due to Tristan Bice~\cite{TristanPersonal}, shows that WAP does not imply AP even in the class of left-cancellative monoids.

\begin{ex}
	Let $M$ be the monoid of transformations of the integers generated by two mappings $a$, $b$, where $\map a \Zee \Zee$ is defined by $a(n) = n+1$ if $n \goe0$ and $a(n) = n-1$ otherwise; $\map b \Zee \Zee$ is defined by $b(n) = n+1$ if $n>0$ and $b(n) = n-1$ otherwise.
	So the only (although essential) difference between $a$ and $b$ is that $a(0) = 1$, while $b(0) = -1$.
	Every element of $M \setminus \sn{\id \Zee}$ is of the form
	$$f = x_{n-1} \cmp \dots x_1 \cmp x_0,$$
	where $x_i \in \dn a b$ for every $i<n$.
	Note that either $f = a^n$ or $f = b^n$, depending on whether $x_0 = a$ or $x_0 = b$.
	This shows that both $a$ and $b$ are amalgamable, however $M$ fails the amalgamation property, because obviously $a$ and $b$ cannot be amalgamated.
\end{ex}

\subsection{Concrete categories}

We now discuss how the results of this note can be interpreted in concrete categories of models and other structures.

First of all, $\fK$ could be a fixed category of finitely generated models of a fixed first-order language while $\fL$ could be the category of all models representable as unions of countable chains in $\ob{\fK}$.
In both cases it is natural to consider embeddings as arrows, where an \define{embedding}{embedding} is an isomorphism onto its image.
It is clear that conditions (L0)--(L2) are satisfied.
In this setting, our results in Sections \ref{SectObjectsGenrkWeakInj}, \ref{SectSlabeHomogenita}, and in particular Theorem~\ref{ThmSzescJedyn}, are extensions of the classical results of \fra~\cite{Fraisse}.
Specifically, if $\ob{\fK}$ is countable up to isomorphism and all the models in $\ob{\fK}$ are countable, then the joint embedding property together with the weak amalgamation property imply the existence of a unique weakly $\fK$-injective model in $\fL$ that might be called the \emph{generic limit} of $\ob{\fK}$.
Note that the property of being hereditary is ignored here.
The main reason is that the weak AP is stable under taking the hereditary closure.
Recall that the \emph{joint embedding property} is simply the property of being directed with respect to embeddings.
If the models in $\fK$ are uncountable (this may happen if the language is uncountable) then we cannot deduce that $\fK$ is locally countable, and indeed $\fK$ might not be weakly dominated by a countable subcategory.
Summarizing, a class $\Emm$ of countable finitely generated models is called a \define{weak \fra\ class}{weak \fra\ class} if it has the joint embedding property, the weak amalgamation property and is essentially countable, namely, has countably many isomorphic types. Once this happens, it is a weak \fra\ category (with embeddings as arrows). This has already been discussed in the recent work~\cite{KraKub}, also in the context of the Banach-Mazur game.
Our Theorem~\ref{ThmSedmPet} in the special case of models summarizes the main results of~\cite{KraKub}.
Recall that if $\Emm$ is a weak \fra\ class then so is its hereditary closure, while if $\Emm$ has the amalgamation property then its hereditary closure may fail the amalgamation property.

\subsection{Projective weak \fra\ theory}

Following Irwin \& Solecki~\cite{IrwSol}, we say that a class of finite nonempty models $\fK$ is a \define{projective \fra\ class}{projective \fra\ class} if it contains countably many types and satisfies the following two conditions:
\begin{enumerate}[itemsep=0pt]
	\item[(1)] For every $X, Y \in \fK$ there exists $Z \in \fK$ having proper epimorphisms onto $X$ and $Y$.
	\item[(2)] Given proper epimorphisms $\map f X Z$, $\map g Y Z$ with $X,Y,Z \in \fK$, there exist $W \in \fK$ and proper epimorphisms $\map{f'}{W}{X}$, $\map{g'}{W}{Y}$ such that $f \cmp f' = g \cmp g'$.
\end{enumerate}
Here, a mapping $\map f A B$ is a \define{proper epimorphism}{proper epimorphism} if it is a surjective homomorphism and satisfies
$$R^B(y_1, \dots, y_n) \iff \; (\exists\; x_1, \dots, x_n \in A)\; R^A(x_1, \dots, x_n) \text{ and } (\forall\; i \loe n) \; y_i = f(x_i)$$
for every $n$-ary relation $R$ (in the language of the models from $\fK$).
It is clear that declaring arrows between $A, B \in \fK$ to be proper epimorphisms from $B$ onto $A$, we obtain a \fra\ category. It is also clear how to change condition (2) above, in order to obtain the projective \emph{weak} amalgamation property.
Of course, the category $\fL \sups \fK$ should consist of all inverse limits of sequences in $\fK$, treated as compact topological spaces with continuous epimorphisms.
It is easy to check that conditions (L0)--(L2) are fulfilled.

As a very concrete example, we may consider $\fK$ to be the class of all finite nonempty sets with no extra structure. Then $\fL$ should be the class of all compact 0-dimensional metrizable spaces. Obviously, $\fK$ is a \fra\ category and its generic limit is the Cantor set.
A much more interesting example (leading to an intriguing topological object, called the \emph{pseudo-arc}) is contained in~\cite{IrwSol}.

\subsection{Uncountable weak \fra\ theory}

There is nothing surprising in extending the theory of generic objects and weak \fra\ sequences to the uncountable setting, namely, working in a category $\fK$ closed under colimits of sequences of length $< \kappa$, where $\kappa$ is an uncountable regular cardinal.
Under certain circumstances, there exists a (unique up to isomorphism) weak \fra\ sequence of length $\kappa$ leading to a generic object in a larger category.
In fact, it suffices to combine the results of Section~\ref{SectDSGdogugo} above with~\cite[Section 3]{Kub40}. Research in this direction has been recently done by Di Liberti~\cite{DiLi}.

In fact, the work~\cite{DroGoe} by Droste \& G\"obel is the first treatment of model-theoretic \fra\ limits from the category-theoretic perspective. Roughly speaking, the authors of~\cite{DroGoe} work in a category $\fL$ having the property that $\lam$-small objects are co-dense and there are not too many of them.
Here, $\lam$ is an infinite regular cardinal.
Under certain natural conditions, $\fL$ contains a special object which is the \fra\ limit of the subcategory of all $\lam$-small objects.
In our case, $\lam = \omega$, however we do not require that $\fK \subs \fL$ consists of \emph{all} $\omega$-small objects. We actually gave necessary and sufficient conditions for the existence of a $\fK$-generic object, assuming conditions (L0)--(L2) only, which are weaker than those of~\cite{DroGoe}. Our main innovation is the concept of a weakly dominating subcategory.
The results of Droste \& G\"obel can be easily extended to the case where the amalgamation property is replaced by its weak version.

\section{Concluding remarks}

We believe that the theory presented in this note will be applicable in various contexts in several areas of mathematics. In fact, it already happened recently that weak \fra\ sequences were used for showing that a certain concrete $C^*$-algebra (called the \emph{Jiang-Su algebra}) is strongly self-absorbing, see~\cite{Saeed}. Roughly speaking, the key idea was considering a category based on the one whose generic limit is the Jiang-Su algebra, showing that this new category has a weak \fra\ sequence and then concluding that its generic limit must be again the Jiang-Su algebra, just because of uniqueness.

As a conclusion, it seems to us that one of the most important aspects of the \fra\ theory (including its weak variant) is uniqueness of the generic object, even though it is obtained almost for free from the Banach-Mazur game.
Nevertheless, \fra\ theory has already provided simple and short proofs of uniqueness of several mathematical objects, like the Gurarii space~\cite{KGur} and the Poulsen simplex~\cite{KubKwi}. Previous arguments were highly non-elementary. 

\subsection{Further research}

First of all, uncountable variants of the Banach-Mazur game presented here should be considered, aiming at finding new connections between ``the finite'' and ``the uncountable''. The game could possibly be played in a non-linear way, building a tree or a more complicated infinite diagram.

Second, the ``real source'' of the weak amalgamation property, namely, structures leading to generic automorphism---this line of research, originated by Ivanov~\cite{Ivanov} and Kechris \& Rosendal~\cite{KechRos} has not been explored yet. We believe that a pure category-theoretic framework can be relatively easily designed here.
%, obtaining new examples.

\printindex

\end{document}